\documentclass[11pt]{article}

\usepackage{xypic}
\usepackage{amsgen}
\usepackage{amsmath}
\usepackage{amstext}
\usepackage{amsbsy}
\usepackage{amsopn}
\usepackage{amsfonts}
\usepackage{amssymb}
\usepackage{eepic}
\usepackage{graphicx}
\usepackage{epsf}
\usepackage{pstricks}
\xyoption{all}

\def\Box{\square}

\def\mapright#1{\smash{\mathop{\longrightarrow}\limits^{#1}}}
\def\tra#1{\smash{\mathop{\mid\kern
-1pt\joinrel\relbar\joinrel\relbar}\limits^{*}_{#1}}}
\def\longtra#1{\smash{\mathop{\mid\kern
-1pt\joinrel\relbar\joinrel\relbar\joinrel\relbar}\limits^{*}_{#1}}}
\def\vlongtra#1{\smash{\mathop{\mid\kern
-1pt\joinrel\relbar\joinrel\relbar\joinrel\relbar\joinrel\relbar}\limits^{*}_{#1}}}
\def\vvlongtra#1{\smash{\mathop{\mid\kern
-1pt\joinrel\relbar\joinrel\relbar\joinrel\relbar\joinrel\relbar\joinrel\relbar}\limits^{*}_{#1}}}
\def\vvvlongtra#1{\smash{\mathop{\mid\kern
-1pt\joinrel\relbar\joinrel\relbar\joinrel\relbar\joinrel\relbar\joinrel\relbar\joinrel\relbar}\limits^{*}_{#1}}}
\def\etra#1{\smash{\mathop{\mid\kern
-1pt\joinrel\relbar\joinrel\relbar}\limits_{#1}}}

\def\vlongrightarrow{\relbar\joinrel\longrightarrow}

\def\vvlongrightarrow{\relbar\joinrel\vlongrightarrow}
\def\vvvlongrightarrow{\relbar\joinrel\vvlongrightarrow}
\def\vvvvlongrightarrow{\relbar\joinrel\vvvlongrightarrow}
\def\vvvvvlongrightarrow{\relbar\joinrel\vvvvlongrightarrow}

\def\longmapright#1{\smash{\mathop{\vlongrightarrow}\limits^{#1}}}
\def\vlongmapright#1{\smash{\mathop{\vvlongrightarrow}\limits^{#1}}}
\def\vvlongmapright#1{\smash{\mathop{\vvvlongrightarrow}\limits^{#1}}}
\def\vvvlongmapright#1{\smash{\mathop{\vvvvlongrightarrow}\limits^{#1}}}
\def\vvvvlongmapright#1{\smash{\mathop{\vvvvvlongrightarrow}\limits^{#1}}}

\def\A{{\cal{A}}}
\def\iff{\Leftrightarrow}
\def\Rw{\Rightarrow}
\def\oo{\overline}

\def\wt{\widetilde}
\def\wh{\widehat}
\def\B{{\cal{B}}}

\def\N{\mathbb{N}}

\def\fix{\mbox{Fix}\,}
\def\irr{\mbox{Irr}\,}

\def\reg{\mbox{Reg}\,}
\def\sing{\mbox{Sing}\,}

\def\ker{\mbox{Ker}}
\def\pref{\mbox{Pref}\,}

\def\max{\mbox{max}}
\def\min{\mbox{min}}
\def\geo{\mbox{Geo}}

\def\TR{{\cal{T}}}
\def\R{{\cal{R}}}

\def\Z{\mathbb{Z}}

\def\p{\varphi}

\def\inv{^{-1}}

\def\bi{\begin{itemize}}
\def\ei{\end{itemize}}
\def\beq{\begin{equation}}
\def\eeq{\end{equation}}

\newtheorem{T}{Theorem}[section]
\newcommand{\bt}{\begin{T}}
\newcommand{\et}{\end{T}}
\newcommand{\ftd}{$\square$\end{T}}

\newtheorem{Proposition}[T]{Proposition}
\newcommand{\bp}{\begin{Proposition}}
\newcommand{\ep}{\end{Proposition}}
\newcommand{\fpd}{$\square$\end{Proposition}}

\newtheorem{Lemma}[T]{Lemma}
\newcommand{\bl}{\begin{Lemma}}
\newcommand{\el}{\end{Lemma}}
\newcommand{\fld}{$\square$\end{Lemma}}

\newtheorem{Corol}[T]{Corollary}
\newcommand{\bc}{\begin{Corol}}
\newcommand{\ec}{\end{Corol}}
\newcommand{\fcd}{$\square$\end{Corol}}

\newtheorem{Result}[T]{Result}
\newcommand{\br}{\begin{Result}}
\newcommand{\er}{\end{Result}}
\newcommand{\frd}{$\square$\end{Result}}

\newtheorem{Example}[T]{Example}
\newcommand{\be}{\begin{Example}}
\newcommand{\ee}{\end{Example}}

\newtheorem{Problem}[T]{Problem}
\newcommand{\bq}{\begin{Problem}}
\newcommand{\eq}{\end{Problem}}

\newcommand{\proof}
   {\par\medbreak\noindent{\bf Proof}.\enspace}

\newcommand{\qed}{%\hfill
$\Box$
\par\bigbreak}

\normalbaselines
\def\abstract#1{\par\bigskip
\begingroup\small
\baselineskip=12truept
\begin{center}ABSTRACT\end{center}
\par\medskip\par\noindent
\null\hfill\hbox{\vbox{\hsize=5truein\noindent#1}}
\hfill\null\par\endgroup\par}

\title{Fixed points of endomorphisms of virtually free groups}
\author{{\bf Pedro V. Silva}\\ $ $\\ {\em Centro de
Matem\'{a}tica, Faculdade de Ci\^{e}ncias, Universidade do
Porto,}\\ {\em R. Campo Alegre 687, 4169-007 Porto, Portugal}\\
{\em e-mail:} pvsilva@fc.up.pt}
\date{\today}

\begin{document}
\maketitle

\begin{center}\small

\bigskip

\end{center}

\abstract{A fixed point theorem is proved for inverse transducers,
  leading to an automata-theoretic proof of the fixed point
subgroup of an endomorphism of a finitely
  generated virtually free group being finitely
  generated. If the endomorphism is uniformly continuous for the
  hyperbolic metric, it is proved that the set of regular fixed points
in the hyperbolic boundary has finitely many orbits under the action of
the finite fixed points. In the automorphism case, it is shown that
these regular fixed points are either exponentially stable attractors
or exponentially stable repellers.}

\section{Introduction}

Throughout the paper, the ambient groups are assumed to be finitely
generated. 

Gersten proved in the 
eighties that the fixed point subgroup of a free group automorphism
$\p$ is
finitely generated \cite{Ger}.
Using a different approach, Cooper gave an alternative
proof, proving also that the fixed points of the
continuous extension of $\p$ to the boundary of the free group is in
some sense finitely generated \cite{Coo}. Bestvina and Handel achieved
in 1992 a 
major breakthrough through their innovative train track techniques,
bounding the rank of the fixed point subgroup and the generating set
for the infinite fixed points \cite{BH}. Their approach was pursued by
Maslakova 
in 2003 to prove that the fixed point subgroup can be effectively
computed \cite{Mas}.

Gersten's result was generalized to further classes of groups and
endomorphisms in subsequent years. Goldstein and Turner extended it to
monomorphisms of free groups \cite{GT}, and later to
arbitrary endomorphisms 
\cite{GT2}. Collins and Turner extended it to automorphisms of free
products of freely indecomposable groups \cite{CT} (see the survey by
Ventura \cite{Ven}). With respect to automorphisms, the widest
generalization is to hyperbolic groups and is due to Paulin \cite{Pau}.

In 2002, Sykiotis extended Collins and Turner's result
to arbitrary endomorphisms of virtually free groups using symmetric
endomorphisms \cite{Syk2} (see
also \cite{Syk} for further results on symmetric endomorphisms). 
In \cite{Sil3}, the author 
generalized Goldstein and Turner's authoma-theoretic
proof to arbitrary endomorphisms of
free products of cyclic groups. In the present paper, this
result is extended to arbitrary endomorphisms of virtually
free groups, providing an automata-theoretic alternative to Sykiotis' result.

This is done by reducing the problem to the rationality of some
languages associated to a finite inverse transducer, and
subsequent application of Anisimov and Seifert's 
Theorem.

Infinite fixed points of automorphisms of free groups were also
discussed by Bestvina and Handel in \cite{BH}. Gaboriau, Jaeger,
Levitt and Lustig remarked in \cite{GJLL} that some of the results on
infinite fixed points would hold for virtually free groups with some
adaptations. 

In \cite{Sil2}, we discussed infinite fixed points for
monomorphisms of free products of cyclic groups, the group case of a
more general setting based on the concept of special confluent
rewriting system. These results are now 
extended to endomorphisms with finite kernel of virtually free groups
(which are precisely the uniformly continuous endomorphisms for the
hyperbolic metric), and we
discuss the dynamical nature of the regular fixed points in the
automorphism case, generalizing the results of \cite{GJLL} on free
groups.

The paper is organized as follows. Section 2 is devoted to
preliminaries on groups and automata. We discuss inverse transducers
in Section 3, proving a useful fixed point theorem. In Section 4 we
prove that the fixed point subgroup is finitely generated for
arbitrary endomorphisms of a (finitely generated)
virtually free group $G$. 

In Section 5 we get a rewriting system with
good properties to represent the elements of $G$, and use it in
Section 6 to construct a simple model for the hyperbolic boundary of
$G$. We study uniformly continuous endomorphisms in Section 7
and prove in Section 8 that the infinite fixed points of such
endomorphisms are in some 
sense finitely generated. 

The classification of the infinite
fixed points of automorphisms is performed in Section 9, and the final
Section 10 includes an example and some open problems.

\section{Preliminaries}

Throughout the whole paper, we assume alphabets to be {\em finite}. 

We start with some group-theoretic definitions.
Given an alphabet $A$, we denote by $A^{-1}$ a set of \emph{formal
  inverses} of $A$,
and write $\widetilde{A} = A\cup A^{-1}$. We extend the mapping $a
\mapsto a\inv$ to an involution of the free monoid $\widetilde{A}^*$
in the obvious way.
As usual, the  \emph{free group
  on} $A$ is the quotient of $\widetilde{A}^*$ by the
congruence generated by the relation 
%\beq
%\label{freegr}
$\{(aa^{-1},1)\mid a \in \widetilde{A}\}.$
%\eeq
We denote by $\theta: \widetilde{A}^* \to F_A$ the canonical
morphism.

Let 
$$R_A = \wt{A}^* \setminus (\cup_{a \in \wt{A}} \; \wt{A}^*aa\inv
  \wt{A}^*)$$
be the subset of all {\em reduced words} in $\wt{A}^*$. It is well
known that, for every $g \in F_A$, $g\theta\inv$ contains a unique
reduced word, denoted by $\oo{g}$. We write also $\oo{u} =
\oo{u\theta}$ for every $u \in \wt{A}^*$. Note that 
the equivalence $u\theta = v\theta \; \iff \; \oo{u} = \oo{v}$
holds for all $u,v \in \wt{A}^*$.

A group $G$ is {\em virtually free} if $G$ has a free subgroup $F$ of
finite index. In view of Nielsen's Theorem, it is well known that $F$
can be assumed to be normal, and is finitely generated if $G$ is
finitely generated itself. Therefore every finitely generated
virtually free group $G$ admits a decomposition as a disjoint union
$$G = F \cup Fb_1 \cup \ldots \cup Fb_m,$$
where $F \unlhd G$ is a free group of finite rank and $b_1,\ldots,b_m
\in G$. 

% \section{Automata}
We shall need also some basic concepts from automata theory:

Let $A$ be a (finite) alphabet. A subset of $A^*$ is called an
$A$-{\em language}. We say that $\A = (Q,q_0,T,\delta)$ is a
(finite) {\em deterministic} $A$-{\em automaton} if: 
\begin{itemize}
\item
$Q$ is a (finite) set;
\item
$q_0 \in Q$ and $T \subseteq Q$;
\item
$\delta: Q \times A \to Q$ is a partial mapping.
\end{itemize}
We extend $\delta$ to a partial mapping $Q \times A^* \to Q$ by
induction through
$$(q,1)\delta = q,\quad (q,ua)\delta = ((q,u)\delta,a)\delta\quad (u
\in A^*,a \in A).$$

When the automaton is clear from the context, we write $qu =
(q,u)\delta$. We can view $\A$ as a directed graph with edges labelled
by letters $a \in A$ by identifying $(p,a)\delta = q$ with the edge
$p \mapright{a} q$. The set of all such edges is denoted by $E(\A)
\subseteq Q \times A \times Q$.

A {\em finite nontrivial path} in $\A$ is a sequence
$$p_0 \mapright{a_1} p_1 \mapright{a_2} \ldots \mapright{a_n} p_n$$
with $(p_{i-1},a_i,p_i) \in E(\A)$ for $i = 1,\ldots,n$. Its {\em label}
is the word $a_1\ldots a_n \in A^*$. It is said to be a {\em
  successful} path if $p_0 = q_0$ and $p_n \in T$. We consider also
the {\em trivial path} $p \mapright{1} p$ for $p \in Q$. It is
successful if $p = q_0 \in T$. 

The {\em language} $L(\A)$ {\em
  recognized by} $\A$
is the set of all labels of successful paths in $\A$. Equivalently, 
$L(\A) = \{ u \in A^* \mid q_0u \in T\}$.
If $(p_{i-1},a_i,p_i) \in E(\A)$ for every $i \in \N$, we may consider also 
the {\em infinite path} 
$$p_0 \mapright{a_1} p_1 \mapright{a_2} p_2 \mapright{a_3} \ldots$$
Its label is the (right) infinite word $a_1a_2a_3\ldots$ We denote by
$A^{\omega}$ the set of all (right) infinite words on the alphabet
$A$, and write also $A^{\infty} = A^* \cup A^{\omega}$. We denote by
$L_{\omega}(\A)$ the set of labels of all infinite 
paths $q_0 \mapright{} \ldots$ in $\A$. 

Given $u \in A^*$ and $\alpha \in A^{\infty}$, we say that $u$ is a
{\em prefix} of $\alpha$ and write $u \leq \alpha$ if $\alpha =
u\beta$ for some $\beta \in A^{\infty}$. By convention, this includes
the case $\alpha \leq \alpha$ for $\alpha \in A^{\omega}$. For every
$n \in \N$, we denote
by $\alpha^{[n]}$ the prefix of length $n$ 
of $\alpha$, applying the
convention that $\alpha^{[n]} = \alpha$ if $n > |\alpha|$. 

It is immediate that $(A^{\infty},\leq)$ is a complete
$\wedge$-semilattice: given $\alpha,\beta \in A^{\infty}$, 
$\alpha\wedge\beta$ is the longest common prefix of $\alpha$ and
$\beta$ (or $\alpha$ if $\alpha = \beta \in A^{\omega}$).
The operator $\wedge$ will play a crucial
role in later sections of the paper. 

The
{\em star} operator on $A$-languages is defined by
$$L^* = \bigcup_{n\geq 0} L^n,$$
where $L^0 = \{ 1 \}$.
An $A$-language $L$ is said to be {\em rational} if $L$ can be obtained
from finite $A$-languages using finitely many times the operators union,
product and star (this is called a {\em rational
  expression}). Alternatively, by Kleene's Theorem \cite[Section III]{Ber},
$L$ is rational if and only if it is 
recognized by a finite deterministic
$A$-automaton $\A$. The
definition through rational expressions generalizes to subsets of an
arbitrary group in the obvious 
way. Moreover, if we fix a homomorphism $\pi:A^* \to
G$, the rational subsets of $G$ are the images by $\pi$ of the
rational $A$-languages. For obvious reasons, we shall be dealing
mostly with matched homomorphisms. A
homomorphism $\pi:\wt{A}^* \to G$ is said to be {\em matched} if
$a\inv\pi = (a\pi)\inv$ for every $a \in A$.
For details on rational languages and subsets, the reader is 
referred to \cite{Ber, Sak}.

We shall need also the following classical result of Anisimov and Seifert:

\bp
\label{anisi}
{\rm \cite[Prop. II.6.2]{Sak}}
Let $H$ be  a subgroup of a group $G$. Then $H$ is a rational subset
of $G$ if and only 
if $H$ is finitely generated.
\ep

We end this section with an elementary observation that will help us to
establish that fixed point subgroups are finitely generated. 

\bp
\label{techno}
Let $\pi:\wt{A}^* \to G$ be a matched epimorphism and let $X \subseteq
G$. Let $\A$ be a finite $\wt{A}$-automaton such that: 
\bi
\item[(i)] $L(\A) \subseteq X\pi\inv$;
\item[(ii)] $L(\A) \cap x\pi\inv \neq \emptyset$ for every $x \in X$.
\ei
Then $X$ is a rational subset of $G$.
\ep

\proof
It follows immediately that $X = (L(\A))\pi$, hence $X$ is a rational
subset of $G$. 
\qed

\section{Inverse transducers}

Given a finite alphabet $A$, we say that $\TR = (Q,q_0,\delta,\lambda)$ is a
(finite) {\em deterministic} $A$-{\em transducer} if: 
\begin{itemize}
\item
$Q$ is a (finite) set;
\item
$q_0 \in Q$;
\item
$\delta: Q \times A \to Q$ and $\lambda: Q \times A \to A^*$ are mappings.
\end{itemize}

As in the automaton case, we may extend $\delta$ to a mapping
$Q \times A^* \to Q$. 
Similarly, we extend $\lambda$ to a mapping $Q \times A^* \to A^*$
through
$$(q,1)\lambda = 1,\quad (q,ua)\lambda =
(q,u)\lambda((q,u)\delta,a)\lambda\quad (u \in A^*,a \in A).$$

When the transducer is clear from the context, we write $qa =
(q,a)\delta$. We can view $\TR$ as a directed graph with edges labelled
by elements of $A\times A^*$ (represented in the form $a|w$) by
identifying $(p,a)\delta = q$, 
$(p,a)\lambda = w$ with the edge
$p \mapright{a|w} q$. The set of all such edges is denoted by $E(\TR)
\subseteq Q \times A \times A^* \times Q$. If $pu = q$ and
$(p,u)\lambda = v$, we write also $p \mapright{u|v} q$ and call it a
path in $\TR$.

It is immediate that, given $u \in A^*$,
there exists exactly one path in $\TR$ of the form 
$q_0 \mapright{u|v} q.$ We write $u\wh{\TR} =
v$, defining thus a mapping $\wh{\TR}: A^* \to A^*$. 

Assume now that $\TR = (Q,q_0,T,\delta,\lambda)$ is a deterministic
$\wt{A}$-transducer such that
$$\mbox{$p \mapright{a|u} q$ is an edge of $\TR$ if and only if 
$q \vvlongmapright{a\inv |u\inv} p$ is an edge of $\TR$.}$$
Then $\TR$ is said to be {\em inverse}.

\bp
\label{induce}
Let $\TR = (Q,q_0,\delta,\lambda)$ be an inverse $\wt{A}$-transducer. Then:
\bi
\item[(i)]
$\delta:Q \times \wt{A}^* \to Q$ induces a mapping $\wt{\delta}:Q
\times F_A \to Q$ by $(q,u\theta)\wt{\delta} = (q,u)\delta$;
\item[(ii)]
$\wh{\TR}: \wt{A}^* \to \wt{A}^*$ induces a partial mapping $\wt{\TR}: F_A \to
F_A$ by $u\theta\wt{\TR} = u\wh{\TR}\theta$. 
\ei
\ep

\proof
(i) Since the free group congruence $\sim$ is generated by the pairs
$(aa\inv,1)$, it suffices to show that $(q,vaa\inv w)\delta =
(q,vw)\delta$ for all $q \in Q$; $v,w \in \wt{A}^*$ and $a \in 
\wt{A}$. 

Since $\delta$ is a full mapping, we have a path 
\beq
\label{induce1}
q \mapright{v|v'} q_1 \mapright{a|u} q_2 \vlongmapright{a\inv|u'} q_3
\mapright{w|w'} q_4
\eeq
in $\TR$. Since $\TR$ is inverse (in particular deterministic), we
must have $u' = u\inv$ and $q_3 = q_1$, hence we also have a path
$$q \mapright{v|v'} q_1 
\mapright{w|w'} q_4$$
and so $(q,vaa\inv w)\delta = q_4 = (q,vw)\delta$ as required.

(ii) Similarly to part (i), it suffices to show that $(vaa\inv
w)\wh{\TR}\theta = 
(vw)\wh{\TR}\theta$ for all $v,w \in \wt{A}^*$ and $a \in
\wt{A}$. 

We consider the path (\ref{induce1}) for $q = q_0$.
Since $u' = u\inv$ and $q_3 = q_1$, we get $$(vaa\inv w)\wh{\TR}\theta
= (v'uu\inv w')\theta = (v'w')\theta = 
(vw)\wh{\TR}\theta$$ as required.
\qed

We prove now one of our main results, generalizing Goldstein and
Turner's proof \cite{GT2} to mappings induced by inverse transducers.

\bt
\label{fixtrans}
Let $\TR$ be a finite inverse $\wt{A}$-transducer and
let $z \in F_A$. Then
$$L = \{ g \in F_A \mid g\wt{\TR} = gz\}$$ 
is rational.
\et

\proof
Write $\TR = (Q,q_0,\delta,\lambda)$.
For every $g \in F_A$, let $P_1(g) = g\inv (g\wt{\TR}) \in F_A$ and
write $q_0g  = (q_0,g)\wt{\delta}$, $P(g) = (P_1(g),q_0g)$. 
Note that $g \in L$ if and only if $P_1(g) = z$.  
We define a deterministic $\wt{A}$-automaton $\A_{\p} = (P,(1,q_0),S,E)$ by
\bi
\item[]
$P = \{ P(g)\mid g \in F_A\}$;
\item[]
$S = P \cap (\{ z \} \times Q)$;
\item[]
$E = \{ (P(g),a,P(ga))\mid g \in F_A,\; a\in \wt{A} \}$.
\ei
Clearly, $\A_{\p}$ is a possibly infinite % complete inverse
automaton. % with a basepoint. 
Note that, since $\TR$ is inverse, we have $qaa\inv = q$ for all $q
\in Q$ and $a \in \wt{A}$. It follows that, whenever
$(p,a,p') \in E$, then also $(p',a\inv,p) \in E$. We say that such edges
are the {\em inverse} of each other. 
% This property, together with
% determinism and trimness (every vertex $P(g)$ lies in some successful
% path, namely $P(1) \mapright{\oo{g}} P(g) \mapright{\oo{g}\inv}
% P(1)$), makes $\A_{\p}$ an {\em inverse} automaton.

Since every $w \in \wt{A}^*$ labels a
unique path $P(1) 
\mapright{w} P(w\theta)$, it follows that
%\beq
%\label{gtholds10}
$$L(\A_{\p}) = L\theta\inv.$$
% \eeq
In view of Proposition \ref{techno}, to prove that $L$ is rational it
suffices to construct a finite subautomaton $\B_{\p}$ of $\A_{\p}$
such that $\oo{L} \subseteq L(\B_{\p})$.
% \beq
% \label{nec1}
% L(\B_{\p}) \cap g\theta\inv \neq \emptyset \mbox{ for every }g \in
% L.
% \eeq

We fix now
$$M = \max\{ |(q,a)\lambda|: q\in Q,\; a \in \wt{A} \},\quad N =
\max\{ 2M+1,|z|\}$$
and 
$$P' = \{ P(g) \in P : |P_1(g)| \leq N \}.$$
Since $A$ is finite, so is $P'$. Given $g \in F_A$, write $g\iota =
\oo{g}^{[1]}$. 
% assume that
% $\oo{g} = a_1\ldots a_n$ with $a_i \in \wt{A}$. We write 
% $$g\iota = \left\{
% \begin{array}{ll}
% a_1&\mbox{ if }n > 0\\
% 1&\mbox{ if }n = 0
% \end{array}
% \right.$$
Given $p = (g,q) \in P$, we write also $p\iota = g\iota$. 
We say that an edge $(p_1,a,p_2) \in E$ is:
\bi
\item
{\em central} if $p_1,p_2 \in P'$;
\item
{\em compatible} if it is not central and $p_1\iota = a$.
\ei

We collect in the following lemma some elementary properties involving
these concepts:

\bl
\label{edges}
\bi
\item[(i)]
There are only finitely many central edges in $\A_{\p}$.
\item[(ii)]
% If $((g_1,q_1),a,(g_2,q_2)) \in E$ is not central, then $g_1\iota =
% (g_1(a\p))\iota$. 
% \item[(iii)]
If $(p_1,a,p_2) \in E$ is not central, then either $(p_1,a,p_2)$ or
$(p_2,a\inv,p_1)$ is compatible.
\item[(iii)]
For every $p \in P$, there is at most one compatible edge leaving $p$.
\ei
\el

\proof
(i) Since $A$ and $P'$ are both finite.

(ii) Assume that $(p_1,a,p_2)$ is neither central nor compatible. Write
$p_1 = (g_1,q_1)$ and $p_2 = (g_2,q_2)$. Suppose that $g_1 = 1$. Then
$g_2 = P_1(a) = a\inv (a\wt{\TR})$ and so $|g_2| \leq 1 + M \leq
N$, in contradiction with $(p_1,a,p_2)$ being non central. 

Thus $\oo{g_1} = bu$ for some
$b \in \wt{A}\setminus \{ a \}$ and $u \in R_A$. On the other hand,
$g_2 = a\inv g_1 
(q_1,a)\lambda$ and so $\oo{g_2} = \oo{a\inv bu
  (q_1,a)\lambda}$. If $|u| < M$, then $|g_1|, |g_2| \leq
2M+1 \leq N$ and $(p_1,a,p_2)$ would be central, a contradiction. Thus
$|u| \geq M \geq |(q_1,a)\lambda|$ and so $g_2\iota = a\inv$. Thus
$(p_2,a\inv,p_1)$ is compatible. 

(iii) Since any compatible edge leaving $p$ must be labelled by
$p\iota$, and $\A_{\p}$ is deterministic.
\qed

A (possibly infinite) path 
$q_0 \mapright{a_1} q_1 \mapright{a_2} \ldots$ in $\A_{\p}$ is:
\bi
\item
{\em central} if all the vertices in 
it are in $P'$;
\item
{\em compatible} if all the edges in 
it are compatible and no intermediate vertex is in $P'$.
\ei

\bl
\label{factpath} 
Let $u \in \oo{L}$. Then there exists a path
$$(1,q_0) = p'_0 \mapright{u_0} p''_0 \mapright{v_1} p_1 \mapright{w_1\inv} p'_1
\mapright{u_1} \ldots \mapright{v_n} p_n \mapright{w_n\inv} p'_n
\mapright{u_n} p''_n \in S$$ 
in $\A_{\p}$ such that:
\bi
\item[(i)]
$u = u_0v_1w_1\inv u_1 \ldots v_nw_n\inv u_n$;
\item[(ii)]
the paths $p'_j
\mapright{u_j} p''_j$ are central;
\item[(iii)]
the paths $p''_{j-1} \mapright{v_j} p_j$ and $p'_j \mapright{w_j} p_j$
are compatible;
\item[(iv)]  
$p_j \notin P'$ if both $v_j$ and $w_j$ are nonempty.
\ei
\el

\proof
Since $S \subseteq P'$ by definition of $N$,  
there exists a path
\beq
\label{factpath1} 
(1,q_0) = p'_0 \mapright{u_0} p''_0 \mapright{x_1} p'_1
\mapright{u_1} \ldots \mapright{x_n} p'_n
\mapright{u_n} p''_n \in S
\eeq
in $\A_{\p}$ such that
$u = u_0x_1u_1 \ldots x_n u_n$ and
the paths $p'_j
\mapright{u_j} p''_j$ (which may be trivial) collect all the occurrences
of vertices in $P'$ (and are therefore central).

By Lemma \ref{edges}(ii), if $(p,a,r)$ occurs in a path $p''_{j-1}
\mapright{x_j} p'_j$, then either $(p,a,r)$ or $(r,a\inv,p)$ is
compatible. On the other hand, since $x_j$ is reduced, it follows
from Lemma \ref{edges}(iii) that $p''_{j-1}
\mapright{x_j} p'_j$ can be factored as  
$$p''_{j-1} \mapright{v_j} p_j \mapright{w_j\inv} p'_j$$
with $p''_{j-1} \mapright{v_j} p_j$ and $p'_j \mapright{w_j} p_j$
compatible. Clearly, (iv) holds since no intermediate vertex of $p''_{j-1}
\mapright{x_j} p'_j$ belongs to $P'$ by construction.
\qed

We say that a compatible path is {\em maximal} if it is infinite or cannot be
extended (to the right) to produce another compatible path.

\bl
\label{adeq}
For every $p \in P'$, there exists in $\A_{\p}$ a unique maximal compatible path
$M_p$ starting at $p$.
\el

\proof
Clearly, every compatible path can be extended to a maximal compatible
path. Uniqueness follows from Lemma \ref{edges}(iii).
\qed 

We define now
$$P'_1 = \{ p \in P' \mid M_p \mbox{ has finitely many distinct edges
} \}$$
and $P'_2 = P' \setminus P'_1$. 
Hence $M_p$ contains no cycles if $p \in P'_2$. 
By Lemma \ref{adeq}, if $M_p$ and $M_{p'}$ intersect at
vertex $r_{pp'}$, then they coincide from $r_{pp'}$ onwards. In particular, if
$M_p$ and $M_{p'}$ intersect, then $p \in P'_1$ if and only if $p' \in
P'_1$. Let 
$$Y = \{ (p,p') \in P'_2 \times P'_2 \mid \mbox{ $M_p$
  intersects }M_{p'}\}.$$
For every $(p,p') \in Y$, let $M_p\setminus M_{p'}$ denote the
(finite) subpath
$p \mapright{} r_{pp'}$ of $M_p$. In
particular, if $p' = p$, $M_p\setminus M_{p'}$ is the trivial path at $p$.

Let $\B_{\p}$ be the subautomaton of $\A_{\p}$ containing:
\bi
\item
all vertices in $P'$ and all central edges;
\item
all edges in the paths $M_p$ $(p \in P'_1)$ and their inverses;
\item
all edges in the paths $M_p\setminus M_{p'}$ $((p,p') \in Y)$ and their
inverses.
\ei
It follows easily from Lemma \ref{edges}(i) and the definitions of $P'_1$ and
$M_p\setminus M_{p'}$ that $\B_{\p}$ is a finite subautomaton of
$\A_{\p}$. As remarked before, it suffices to show that $\oo{L}
\subseteq L(\B_{\p})$.

Let $u \in \oo{L}$. Since $\B_{\p}$ contains all the central edges of
$\A_{\p}$, it suffices to show that all subpaths $$p''_{j-1} \mapright{v_j}
p_j \mapright{w_j\inv} p'_j$$ appearing in the factorization provided by
Lemma \ref{factpath} are paths in $\B_{\p}$. 

Without loss of generality, we may assume that $v_j \neq 1$. If $w_j =
1$, then $p''_{j-1} \in P'_1$ and we are done, hence we may assume
that also $w_j \neq 1$. Now, if
one of the vertices
$p''_{j-1}, p'_j$ is in $P'_1$, so is the other and we are done since
$\B_{\p}$ contains all the edges in the paths $M_p$ $(p \in P'_1)$ and
their inverses. Hence we may assume that $p''_{j-1}, p'_j \in
P'_2$. It follows that $p_j =
r_{p''_{j-1},p'_j}$ (since 
$v_jw_j\inv \in R_A$, the paths $M_{p''_{j-1}}$ and $M_{p'_j}$ cannot meet before
$p_j$). Thus $p''_{j-1} \mapright{v_j}
p_j$ is $M_{p''_{j-1}} \setminus M_{p'_j}$ and $p'_{j} \mapright{w_j}
p_j$ is $M_{p'_j} \setminus M_{p''_{j-1}}$, and so these are also
paths in $\B_{\p}$ as required.
\qed

\section{The fixed point subgroup}

We can now produce an automata-theoretic proof to Sykiotis' theorem:

\bt
\label{gtholds}
{\rm \cite[Proposition 3.4]{Syk2}}
Let $\p$ be an endomorphism of a finitely generated
virtually free group. Then {\rm Fix}$\,\p$ is finitely
generated. 
\et

\proof
We consider a decomposition of $G$ as a disjoint union
\beq
\label{gtholds1}
G = Fb_0 \cup Fb_1 \cup \ldots \cup Fb_m,
\eeq
where $F = F_A \unlhd G$ is a free group with $A$ finite and $b_0,\ldots,b_m
\in G$ with $b_0 = 1$. 

% We fix an extended alphabet $B = A \cup \{ b_1,\ldots,b_m\}$ and let
% $\pi: \wt{B}^* \to G$ denote the canonical epimorphism. 

Let $\p_0:F_A \to F_A$ and $\eta:F_A \to \{ 0,\ldots,m\}$ be defined
by
$$g\p = (g\p_0)b_{g\eta} \quad (g \in F_A).$$
Since the decomposition (\ref{gtholds1}) is disjoint, $g\p_0$ and
$g\eta$ are both uniquely determined by $g\p$, and so both mappings
are well defined. 

Write $Q = \{ 0,\ldots,m\}$. For all $i \in Q$ and $a\in \wt{A}$, we
have $b_i(a\p) = h_{i,a} b_{(i,a)\delta}$ for some (unique)
$h_{i,a} \in F_A$ and $(i,a)\delta \in Q$.
It follows that, for every $j \in Q$,  $\A_j =
(Q,0,j,\delta)$ is a well-defined finite deterministic
$\wt{A}$-automaton. 
We define also a finite deterministic
$\wt{A}$-transducer $\TR = (Q,0,\delta,\lambda)$ by taking
$(i,a)\lambda = \oo{h_{i,a}}$ for all $i \in Q$ and $a\in \wt{A}$. 

Assume that $$i \vlongmapright{a|\oo{h_{i,a}}} (i,a)\delta = j$$ is an
edge of $\TR$. Then $b_i(a\p) = h_{i,a}b_j$ and so also
$$b_i = b_i(a\p)(a\inv\p) = h_{i,a}b_j(a\inv\p) =
h_{i,a}h_{j,a\inv}b_{(j,a\inv)\delta}.$$
This yields $h_{i,a}h_{j,a\inv} = 1$ and $(j,a\inv)\delta = i$, thus
there is an edge $j \vvvlongmapright{a\inv|\oo{h_{i,a}}\inv}
(j,a\inv)\delta = i$ in $\TR$ and so $\TR$ is an inverse
transducer. We claim that $\wt{\TR} = \p_0$. Indeed, let $g =
a_1\ldots a_n$ $(a_i \in \wt{A_i})$. Then there exists a (unique) path
in $\TR$ of the form
$$0 = i_0 \vvvlongmapright{a_1|\oo{h_{i_0,a_1}}} i_1
\vvvlongmapright{a_2|\oo{h_{i_1,a_2}}} \ldots
\vvvvlongmapright{a_n|\oo{h_{i_{n-1},a_n}}} i_n.$$ 
Moreover, $i_j = (i_{j-1},a_j)\delta$ for $j = 1,\ldots,n$. 
It follows that
$$\begin{array}{lll}
g\p&=&b_{i_0}(a_1\p)\ldots(a_n\p) =
h_{i_0,a_1}b_{i_1}(a_2\p)\ldots(a_n\p) =
h_{i_0,a_1}h_{i_1,a_2}b_{i_2}(a_3\p)\ldots(a_n\p)\\ &=&\ldots =
h_{i_0,a_1} \ldots h_{i_{n-1},a_n}b_{i_n} 
\end{array}$$
and so
$$g\p_0 = h_{i_0,a_1} \ldots h_{i_{n-1},a_n} = (\oo{h_{i_0,a_1}} \ldots
\oo{h_{i_{n-1},a_n}})\theta = g\wt{\TR}.$$
Thus $\wt{\TR} = \p_0$.

Note that we have also shown that $g\eta = i_n = (0,a_1\ldots
a_n)\delta$, hence  
\beq
\label{gtholds2}
L(\A_j) = \{ u \in \wt{A}^* \mid u\theta\eta = j\}.
\eeq

Next let 
$$Y = \{ (i,j) \in Q \times Q \mid b_j(b_i\p) \in F_Ab_i\}.$$
For every $(i,j) \in Y$, let $z_{i,j} \in F_A$ be such that
$b_j(b_i\p) = z_{i,j}b_i$ and define
$$X_{i,j} = \{ g \in F_A \mid gb_i \in \fix\p \mbox{ and } g\eta =
j\}.$$ 
We claim that $X_{i,j}$ is a rational subset of $F_A$ for every $(i,j)
\in Y$. 
Indeed, $(gb_i)\p = (g\p)(b_i\p) =
(g\p_0)b_{g\eta}(b_i\p)$. Hence
$$\begin{array}{lll}
X_{i,j}&=&\{ g \in F_A \mid (g\p_0)b_{j}(b_i\p) = gb_i \mbox{ and } g\eta =
j\} = \{ g \in F_A \mid (g\p_0)z_{i,j}b_i = gb_i \mbox{ and } g\eta =
j\}\\
&=&\{ g \in F_A \mid g\p_0 = gz_{i,j}\inv \} \cap \{ g \in F_A
\mid g\eta = j\}.
\end{array}$$
Writing
$$L_{i,j} = \{ g \in F_A \mid g\p_0 = gz_{i,j}\inv \},$$
it follows from (\ref{gtholds2}) that
$X_{i,j} = L_{i,j} \cap (L(\A_j))\theta$. Since $\p_0 = \wt{\TR}$, it
follows from Theorem \ref{fixtrans} that $X_{i,j}$ is an intersection
of two rational subsets of $F_A$, hence rational itself (see
\cite[Corollary III.2.10]{Ber}). 

Now it is easy to check that
\beq
\label{gtholds3}
\fix\p = \cup_{i \in Q} \; (\cup\{ X_{i,j} \mid (i,j) \in Y\})b_i.
\eeq
Indeed, for every $(i,j) \in Y$, we have $X_{i,j}b_i \subseteq \fix\p$
by definition of $X_{i,j}$. Conversely, let $gb_i \in \fix\p$ for some
$g\in F_A$ and $i \in Q$. Then
$gb_i = (gb_i)\p = (g\p_0)b_{g\eta}(b_i\p)$ and so
$b_{g\eta}(b_i\p) \in F_Ab_i$. Hence $(i,g\eta) \in Y$. Since $g
\in X_{i,g\eta}$, (\ref{gtholds3}) holds. Since the $X_{i,j}$ are
rational subsets of $F_A$ and therefore of $G$, it follows that
$\fix\p$ is a rational subset of $G$ and therefore finitely generated by
Proposition \ref{anisi}.
\qed

\section{A good rewriting system}
\label{sfive}

We recall that a (finite) {\em rewriting system} on $A$ is a (finite)
subset $\R$ of $A^* \times A^*$. 
Given $u,v \in A^*$, we write $u
\mapright{}_{\R} v$ if there exist $(r,s) \in \R$ and $x,y \in A^*$ such
that $u = xry$ and $v = xsy$. The reflexive and transitive closure of
$\mapright{}_{\R}$ is denoted by $\mapright{}_{\R}^*$.

We say that $\R$ is:
\bi
\item
{\em length-reducing} if $|r| > |s|$ for every $(r,s) \in \R$;
\item
{\em length-nonincreasing} if $|r| \geq |s|$ for every $(r,s) \in \R$;
\item
{\em noetherian} if, for every $u \in A^*$, there is a bound on the
length of a chain $$u \mapright{}_{\R} v_1 \mapright{}_{\R} \ldots
\mapright{}_{\R} v_n;$$
\item
{\em confluent} if, whenever $u \mapright{}_{\R}^* v$ and $u
\mapright{}_{\R}^* w$, there exists some $z \in A^*$ such that $v
\mapright{}_{\R}^* z$ and $w \mapright{}_{\R}^* z$.
\ei
A word $u \in A^*$
is an {\em irreducible} if no $v \in A^*$ satisfies $u \mapright{}_{\R}
v$. We denote by $\irr\R$ the set of all irreducible words in $A^*$
with respect to $\R$.

We introduce now some basic concepts and results from the
theory of hyperbolic groups. For
details on this class of groups, the reader is referred to \cite{GH}.

Let $\pi:\widetilde{A}^* \to G$ be a matched epimorphism with $A$ finite. 
The {\em Cayley graph} $\Gamma_A(G)$ of $G$ with respect to $\pi$ has
vertex set $G$ and edges $(g,a,g(a\pi))$ for 
all $g \in G$ and $a \in \widetilde{A}$. We say that a path $p \mapright{u}
q$ in $\Gamma_A(G)$ is a {\em geodesic} if it has shortest length among all
the paths connecting $p$ to $q$ in $\Gamma_A(G)$. We denote by
$\geo_A(G)$ the set of labels of all geodesics in $\Gamma_A(G)$. Note
that, since $\Gamma_A(G)$ is vertex-transitive, it is irrelevant
whether or not we fix a basepoint.

The {\em geodesic distance} $d_1$ on $G$ is defined by taking $d_1(g,h)$
to be the length of a geodesic from $g$ to $h$. 
Given $X \subseteq G$
nonempty and $g \in G$, we define
$$d_1(g,X) = \min\{ d_1(g,x) \mid x \in X\}.$$
A {\em geodesic triangle} in $\Gamma_A(G)$ is a collection of three geodesics
$$P_1: g_1 \mapright{} g_2,\quad P_2: g_2 \mapright{} g_3, \quad P_3:
g_3 \mapright{} g_1$$ 
connecting three vertices $g_1,g_2,g_3 \in G$. Let $V(P_i)$ denote the set
of vertices occurring in the path $P_i$. We say that $\Gamma_A(G)$ is
$\delta$-{\em hyperbolic} for some $\delta \geq 0$ if
$$\forall g \in V(P_1) \hspace{.7cm} d_1(g,V(P_2)\cup V(P_3)) < \delta$$
holds for every geodesic triangle $\{ P_1,P_2,P_3\}$ in in
$\Gamma_A(G)$. If this happens for some $\delta$, we say that $G$ is {\em
  hyperbolic}. It is well known that the concept is independent from
both alphabet and matched epimorphism, but the hyperbolicity
constant $\delta$ may change. Virtually free groups are among the most
important examples of hyperbolic groups.

We now use a theorem of Gilman, Hermiller, Holt and Rees \cite{GHHR} to
prove the following result:

\bl
\label{bouge}
Let $G$ be a finitely generated virtually free group. Then there
exist a finite alphabet $A$, a matched epimorphism $\pi:\widetilde{A}^*
\to G$ and a positive integer $N_0$ such that, for all $u \in$ {\rm
  Geo}$_A(G)$ and $v \in \widetilde{A}^*$:
\bi
\item[(i)] there exists some $w \in$  {\rm
  Geo}$_A(G)$ such that $w\pi = (uv)\pi$ and $|u\wedge w| \geq
|u|-N_0|v|$;
\item[(ii)] there exists some $z \in$  {\rm
  Geo}$_A(G)$ such that $z\pi = (vu)\pi$ and $|u\inv \wedge z\inv| \geq
|u|-N_0|v|$.
\ei
\el

\proof
(i) By \cite[Theorem 1]{GHHR}, there
exists a finite alphabet $A$, a matched epimorphism $\pi:\widetilde{A}^*
\to G$ and a finite length-reducing rewriting system $\R$ such that
$\geo_A(G) = \irr\R$. The authors also prove that this property
characterizes (finitely generated) virtually free groups.

Let $N_0 = \max\{ |r|: (r,s) \in \R \}$. Suppose that
$$uv = w_0 \mapright{}_{\R} w_1 \mapright{}_{\R} \ldots
\mapright{}_{\R} w_n = w$$
is a sequence of reductions leading to a geodesic $w$. Then
$(wv\inv)\pi = u\pi$ and since $u$ is a geodesic we get $|w| \geq |u|
- |v|$. Since $\R$ is length-reducing, this yields $n \leq |u|-|w| \leq |v|$. 

Trivially, $|u\wedge w_0| \geq |u|$. Since $u\wedge w_{i-1} \in
\geo_A(G)$, it is immediate that $|u\wedge
w_i| > |u\wedge w_{i-1}| - N_0$ and so $$|u\wedge w| = |u\wedge w_{n}|
\geq |u| - nN_0 \geq |u| - N_0|v|.$$

(ii) The inverse of a geodesic is still a geodesic. By applying (i) to
$u\inv$ and $v\inv$, we get $(u\inv v\inv)\pi = x\pi$ for some $x \in
\geo_A(G)$ satisfying $|u\inv \wedge x| \geq
|u\inv|-N_0|v\inv|$. Then we take $z = x\inv$.
\qed

We assume for the remainder of the paper that 
$G$ is a finitely generated virtually free group, $\pi:\widetilde{A}^*
\to G$ a  matched epimorphism and $N_0$ a positive integer
satisfying the conditions of Lemma \ref{bouge}. Since $G$ is
hyperbolic, it follows from \cite[Theorem 3.4.5]{Eps} that $\geo_A(G)$ is an
automatic structure for $G$ with respect to $\pi$ (see \cite{Eps} for
definitions),
and so the {\em 
  fellow traveller property} holds for some constant $K_0 > 0$ (which
can be taken as $2(\delta +1)$, if $\delta$ is the hyperbolicity constant). This
amounts to say that
% \beq
% \label{fellow}
$$\forall u,v \in \geo_A(G) \: (d_1(u\pi,v\pi) \leq 1 \Rw \forall n
\in \N\; 
d_1(u^{[n]}\pi, v^{[n]}\pi) \leq K_0).$$
% \eeq
We fix a total ordering of $\widetilde{A}$. The {\em shortlex
  ordering} of $\widetilde{A}^*$ is defined by
$$u \leq_{sl} v \mbox{ if }\left\{
\begin{array}{c}
|u| < |v|\\
\mbox{ or }\\
|u| = |v| \mbox{ and } u = wau', v = wbv' \mbox{ with } a<b \mbox{ in
}\widetilde{A}
\end{array}
\right.$$
This is a well-known well-ordering of $\widetilde{A}^*$, compatible
with multiplication on the left and on the right.
Let
$$L = \{ u \in \geo_A(G) \mid u \leq_{sl} v \mbox{ for every } v \in
u\pi\pi\inv \}.$$
By \cite[Theorem 2.5.1]{Eps}, $L$ is also an automatic structure for $G$ with
respect to $\pi$. We note that
$L$ is {\em factorial} (a factor of a word in $L$ is still in $L$).

Given $g \in G$, let $\oo{g}$ denote the unique word of
$L$ representing $g$. This corresponds precisely to free group
reduction if $G = F_A$ and $\pi = \theta$. Since we shall not need
free group reduction from now on, we write also $\oo{u} = \oo{u\pi}$
for every $u \in \wt{A}^*$ to simplify notation.

\bt
\label{spn}
Consider the finite rewriting system $\R'$ on $A$ defined by
$$\R' = \{ (u,\oo{u}): u \in \wt{A}^*,\; |u| \leq K_0N_0 + 1, \; u \neq
\oo{u} \}.$$ Then:
\bi
\item[(i)] $\R'$ is length-nonincreasing, noetherian and confluent;
\item[(ii)] $\irr \R' = L$.
\ei
\et

\proof
(i) $\R'$ is trivially length-nonincreasing, and noetherian folows
from 
\beq
\label{spn2}
(u,\oo{u}) \in \R' \Rw u >_{sl} \oo{u}
\eeq
and $\wt{A}^*$ being
well-ordered by $\leq_{sl}$. 

Next we show that
\beq
\label{spn1}
u \mapright{}_{R'}^* \oo{u} \hspace{1cm}\mbox{holds for every } u \in
\wt{A}^*.
\eeq

We use induction on $|u|$. The case $|u| \leq K_0N_0 + 1$ follows from the
definition of $\R'$, hence assume that $|u| > K_0N_0 + 1$ and (\ref{spn1})
holds for shorter words. Write $u = avb$ with $a,b \in
\wt{A}$. If $av \notin L$, we have $u
\mapright{}_{R'}^* \oo{av}b$ and $\oo{u} =
\oo{\oo{av}b}$, hence $u \mapright{}_{R'}^* \oo{u}$ follows from $\oo{av}y
\mapright{}_{R'}^* \oo{\oo{av}y}$. Hence we may assume that $av \in
L$. 

Suppose that $u \notin \geo_A(G)$. By Lemma \ref{bouge}(i), there
exists some $w \in \geo_A(G)$ such that $w\pi = (avb)\pi$ and $|av\wedge w| \geq
|av|-N_0 \geq K_0N_0 + 1 - N_0 > 0$. Hence we may write $w = aw'$ and we
get $(vb)\pi = (a\inv w)\pi = w'\pi$. Since $|w'| < |vb|$ due to $u
\notin \geo_A(G)$, we get $|\oo{vb}| < |vb|$ and so we may apply twice
the induction hypothesis to get 
$$u = avb \mapright{}_{R'}^* a\oo{vb} \mapright{}_{R'}^* \oo{a\oo{vb}}
= \oo{u}.$$
Hence we may assume that $u \in \geo_A(G)$.
We claim that $\oo{u}^{[1]} = a$.
Let $p = K_0N_0 + 1$. Since $u,\oo{u} \in \geo_A(G)$ and $u\pi = \oo{u}\pi$,
the fellow traveller
property yields $d_1(u^{[p]}\pi,\oo{u}^{[p]}\pi) \leq K_0$ and so
$u^{[p]}\pi = (\oo{u}^{[p]}x)\pi$ for some $x$ of length $\leq
K_0$. Thus, by Lemma \ref{bouge}(i), there exists some $w \in
\geo_A(G)$ such that $w\pi = (\oo{u}^{[p]}x)\pi = u^{[p]}\pi$ and
$$|\oo{u}^{[p]}\wedge w| \geq |\oo{u}^{[p]}|-N_0|x| \geq p-K_0N_0 =
1,$$
hence $\oo{u}^{[1]} = w^{[1]}$. Now $av \in L$ by assumption, hence
$u^{[p]} \in L$ and so $u^{[p]} = \oo{u^{[p]}}$. Since $w\pi = u^{[p]}\pi$ and $w \in
\geo_A(G)$, we get $a = u^{[1]} \leq w^{[1]} = \oo{u}^{[1]}$ in 
$(\wt{A},\leq)$. On the other hand, $\oo{u} \leq_{sl} u$ yields
$\oo{u}^{[1]} \leq a$ in $(\wt{A},\leq)$ and so $\oo{u}^{[1]} = a$ as
claimed.

Now it
follows easily that $\oo{u} = a\oo{a\inv u} = a\oo{vb}$ and the
induction hypothesis yields $vb \mapright{}_{R'}^* \oo{vb}$
and therefore $u = avb \mapright{}_{R'}^* a\oo{vb} =
\oo{u}$. Therefore (\ref{spn1}) holds.

Assume now that $u \mapright{}_{R'}^* v$ and $u \mapright{}_{R'}^*
w$. By (\ref{spn1}), we get $v \mapright{}_{R'}^* \oo{v} = \oo{u}$ and $w
\mapright{}_{R'}^* \oo{w} = \oo{u}$, hence $\R'$ is confluent.

(ii) It follows from (\ref{spn1}) that $\irr\R' \subseteq L$. The
converse inclusion follows from the implication
$$u \mapright{}_{R'} v \Rw u >_{sl} v,$$
which follows in turn from (\ref{spn2}).
\qed

We establish now some technical results which will be useful in later
sections:

\bl
\label{fstbound}
Let $u,v \in L$ and let $w \in \widetilde{A}^*$ be such
  that $vw \in$ {\rm 
  Geo}$_A(G)$ and $(vw)\pi = u\pi$. Then $|u\wedge v| \geq |v|
-K_0N_0$.
\el

\proof
Let $k = |v|$ and write $u = u^{[k]}u'$. Since $v = (vw)^{[k]}$, it
follows from the fellow traveller
property that $d_1(v\pi, u^{[k]}\pi) \leq K_0$, hence we may write
$v\pi = (u^{[k]}z)\pi$ with $|z| \leq K_0$. Since $u^{[k]}$ is itself
a geodesic, it follows from Lemma \ref{bouge}(i) that there exists a
geodesic $u^{[p]}z'$ satisfying $(u^{[p]}z')\pi = (u^{[k]}z)\pi =
v\pi$ and $$p = |u^{[k]}\wedge u^{[p]}z'| \geq |u^{[k]}|-N_0|z|
\geq |v|-K_0N_0.$$ 
Now $v \in L$ yields $v \leq_{sl} u^{[p]}z'$ and so $v^{[p]} \leq_{sl}
u^{[p]}$. On the other hand, $u \in L$ yields $u \leq_{sl} vw$ and so
$u^{[p]} \leq_{sl} v^{[p]}$. Thus $u^{[p]} = v^{[p]}$ and so $|u\wedge
v| \geq p \geq |v|
-K_0N_0$.
\qed

\bp
\label{unibo}
\bi
\item[(i)] Let $uv \in L$ and let $w \in \widetilde{A}^*$ be such that $|v|
  \geq K_0N_0 + N_0|w|$. Then $\oo{uvw} = u\oo{vw}$.
\item[(ii)] Let $u \in \widetilde{A}^*$ and let $vw,vw' \in L$. Then
  $|\oo{uvw} \wedge \oo{uvw'}| \geq |v| - K_0N_0 - N_0|u|$.  
\ei
\ep

\proof
(i) Write $v = v_1v_2$ with $|v_2| = N_0|w|$. By Lemma \ref{bouge}(i),
there exists some $uv_1z \in \geo_A(G)$ such that $(uv_1z)\pi =
(uvw)\pi$. Let $x = \oo{uvw}$. By Lemma \ref{fstbound}, we get $|x
\wedge uv_1| \geq |uv_1| - K_0N_0$. Since $|v_1| = |v| - |v_2| \geq
K_0N_0$, then $u \leq x$ and we may write $x = uy$ for some $y$. Since
$L$ is factorial, we have $y \in L$. In view of $y\pi = (u\inv x)\pi =
(vw)\pi$, we get $y = \oo{vw}$ and so $\oo{uvw} = u\oo{vw}$.

(ii) We may assume that $|v| > K_0N_0 + N_0|u|$.
Write $v = v_1v_2$ with $|v_1| = N_0|u|$. Let $x = \oo{uv_1}$ and
write $p = |x|+|v_2|$.
By the proof of Lemma \ref{bouge},
we have $xv_2w, xv_2w' \in \geo_A(G)$. 

Let $y = \oo{uvw}$. Since $(xv_2w)\pi = y\pi$, it
follows from the fellow traveller
property that $d_1((xv_2)\pi, y^{[p]}\pi) \leq K_0$, hence we may write
$(xv_2)\pi = (y^{[p]}s)\pi$ with $|s| \leq K_0$. Since
$y^{[p]}$ is itself 
a geodesic, it follows from Lemma \ref{bouge}(i) that there exists a
geodesic $y^{[p-K_0N_0]}s'$ satisfying $(y^{[p-K_0N_0]}s')\pi = (y^{[p]}s)\pi =
(xv_2)\pi$. To complete the proof, it suffices to show that
\beq
\label{unibo1}
|y \wedge \oo{xv_2}| \geq p-K_0N_0.
\eeq
Indeed, together with the corresponding inequality for $y' = \oo{uvw'}$,
this implies 
$$|\oo{uvw} \wedge \oo{uvw'}| \geq p-K_0N_0 \geq |v_2|-K_0N_0 = |v|
- K_0N_0 - N_0|u|$$ and we obtain the desired inequality.

To prove (\ref{unibo1}), we consider the geodesic
$y^{[p-K_0N_0]}s'$. Since $(y^{[p-K_0N_0]}s')\pi = (xv_2)\pi$, we get
$\oo{xv_2} \leq_{sl} y^{[p-K_0N_0]}s'$ and so
$\oo{xv_2}^{[p-K_0N_0]} \leq_{sl} y^{[p-K_0N_0]}$. On the other hand, $xv_2w$ is
also a geodesic, hence $y = \oo{uvw} = \oo{xv_2w} \leq_{sl} \oo{xv_2}w$
yields $y^{[p-K_0N_0]} \leq_{sl} \oo{xv_2}^{[p-K_0N_0]}$.
Therefore $y^{[p-K_0N_0]} = \oo{xv_2}^{[p-K_0N_0]}$ and so
(\ref{unibo1}) holds as required.
\qed

\section{A new model for the boundary}
\label{ssix}

We can now present a new model for the boundary of a finitely generated
virtually free group which will prove itself fit do study infinite
fixed points in forthcoming sections. The notion of boundary is indeed
one of the important features associated to hyperbolic groups. To
present it, we shall define a second distance 
in $G$ by means of the {\em Gromov product} (taking 1 as
basepoint). We keep all the notation introduced in Section
\ref{sfive}. In particular, $G$ is a finitely generated virtually free
group and $L = \irr\R'$.

Given $g,h \in G$, we 
define 
$$(g|h) = \frac{1}{2}(d_1(1,g) + d_1(1,h)-d_1(g,h)).$$
Fix $\varepsilon > 0$ such that $\varepsilon\delta \leq
\frac{1}{5}$. Write $z = e^{\varepsilon}$ and define 
$$\rho(g,h) =
\left\{
\begin{array}{ll}
z^{-(g|h)}&\mbox{ if } g\neq h\\
0&\mbox{ otherwise}
\end{array}
\right.$$
for all $g,h \in G$. In general, $\rho$ is not a distance
because it fails the triangular inequality. This problem is overcome
by defining
$$d_2(g,h) = \inf\{ \rho(g_0,g_1) + \ldots + \rho(g_{n-1},g_n) \mid
g_0 = g, \; g_n = h;\; g_1,\ldots,g_{n-1} \in G \}.$$
By \cite[Proposition 5.16]{Vai} (see also \cite[Proposition
7.10]{GH}), $d_2$ is a distance on $G$ and the inequalities
\beq
\label{ineq}
\frac{1}{2}\rho(g,h) \leq d_2(g,h) \leq \rho(g,h)
\eeq
hold for all $g,h \in G$. 

In general, the metric space $(G,d_2)$ is not complete. Its completion
$(\wh{G},\wh{d_2})$ is essentially unique, and $\partial G = \wh{G}
\setminus G$ is the {\em boundary} of $G$. The elements of the
boundary admit several standard descriptions, such as equivalence
classes of rays (infinite words whose finite factors are geodesics)
when two rays are equivalent if the Hausdorff distance between them is
finite \cite[Section 7.1]{GH}. We won't need precise definitions for
these concepts or $\wh{d_2}$ since, as we shall see next, we can get
a simpler description of $\wh{G}$ for virtually free groups.

\bl
\label{tribound}
There exists some $M_0 > 0$ such that, for all $g,h \in G$:
\bi
\item[(i)] $|\oo{g}| \leq |\oo{g}\wedge \oo{gh}| + K_0N_0+ N_0|\oo{h}|$;
\item[(ii)] $d_1(g,h) \geq \frac{|\oo{g}|-|\oo{g}\wedge \oo{h}|}{N_0}
  - K_0$;
\item[(iii)] $|\oo{g} \wedge \oo{h}| \leq
(g|h) \leq |\oo{g} \wedge \oo{h}| + M_0$.
\ei
\el

\proof
(i) By applying Lemma \ref{bouge} to the product $\oo{g}\oo{h}$, there
exists some factorization $\oo{g} = vz$ and some geodesic $vw \in
(gh)\pi\inv$ such that $|v| \geq |\oo{g}|-N_0|\oo{h}|$. Now we apply
Lemma \ref{fstbound} to $u = \oo{gh}$ and $vw$ to get
$|u\wedge v| \geq |v|
-K_0N_0$. Hence
$$|\oo{g}\wedge \oo{gh}| = |u\wedge v| \geq |v|
-K_0N_0 \geq |\oo{g}| - N_0|\oo{h}| -K_0N_0$$
and (i) holds.

(ii) Let $u = \oo{g} \wedge \oo{h}$. Applying (i) to $g$ and $g\inv
h$, and in view of $d_1(g,h) = |\oo{g\inv h}|$, we get
$$|\oo{g}| \leq |\oo{g}\wedge \oo{h}| + K_0N_0+ N_0d_1(g,h)$$
and so (ii) holds.

(iii) We define $M_0 = \delta + (2\delta+1+K_0)N_0 -\frac{1}{2}$,
assuming that $\geo_A(G)$ is 
$\delta$-hyperbolic. 
Let $u = \oo{g} \wedge \oo{h}$, and write $\oo{g} = uv$, $\oo{h} =
uw$. 
It is easy to check that
$$(g|h) = \frac{1}{2}(d_1(1,g)+d_1(1,h)-d_1(g,h)) = \frac{1}{2}(|u| +
d_1(u\pi,g)+ |u| + d_1(u\pi,h)-d_1(g,h)).$$ 
Since $d_1(g,h) \leq d_1(g,u\pi)+ d_1(u\pi,h)$, we get $|\oo{g} \wedge
\oo{h}| = |u| \leq (g|h)$. 

Consider now the geodesic triangle determined by the paths
$$P_1:u\pi \mapright{v} g,\quad P_2:u\pi \mapright{w} h,\quad P_3:g
\vlongmapright{\oo{g\inv h}} h.$$ 
Since $\geo_A(G)$ is $\delta$-hyperbolic, then
\beq
\label{tribound1}
d_1(q,V(P_1)\cup V(P_2)) < \delta \; \mbox{
  for every }q \in V(P_3).
\eeq
Assume that $P_3:g = q_0 \mapright{a_1} \ldots \mapright{a_n} q_n = h$
with $a_i \in \widetilde{A}$. Since $d_1(q_0,V(P_1)) = 0 < \delta$ and
$d_1(q_n,V(P_2)) = 0 < \delta$, it follows from (\ref{tribound1}) that
there exist some $j \in \{ 0,\ldots, n-1\}$ and $p_1 \in V(P_1)$,
$p_2 \in V(P_2)$ such that $d_1(q_j,p_1), d_1(q_{j+1},p_2) \leq
\delta$. 
Since $P_1$ and $P_2$ are geodesics, we get
$$\begin{array}{lll}
(g|h)&=&\frac{1}{2}(d_1(1,g)+d_1(1,h)-d_1(g,h))\\
&=&\frac{1}{2}(|u| + d_1(u\pi,p_1)+d_1(p_1,g)\\
&&\hspace{1cm} + |u| + d_1(u\pi,p_2)+d_1(p_2,h)
-d_1(g,q_j)-1-d_1(q_{j+1},h))\\
&=&|\oo{g} \wedge \oo{h}| + \frac{1}{2}(d_1(u\pi,p_1)+d_1(u\pi,p_2))\\
&&\hspace{1cm} + 
\frac{1}{2}(d_1(p_1,g)-d_1(g,q_j)) +
\frac{1}{2}(d_1(p_2,h)-d_1(q_{j+1},h)) -\frac{1}{2}.
\end{array}$$
Since $d_1(p_1,g) \leq d_1(p_1,q_j) + d_1(q_j,g) \leq \delta +
d_1(q_j,g)$, we have
$$\frac{1}{2}(d_1(p_1,g)-d_1(g,q_j)) \leq \frac{\delta}{2}.$$
Similarly, 
$$\frac{1}{2}(d_1(p_2,h)-d_1(q_{j+1},h)) \leq \frac{\delta}{2}.$$
Out of symmetry, it suffices to show that $d_1(u\pi,p_1) \leq (2\delta
+1+K_0)N_0$.

Applying (ii) to $p_1$ and $p_2$, we get 
$$d_1(p_1,p_2) \geq \frac{|\oo{p_1}|-|\oo{p_1}\wedge \oo{p_2}|}{N_0}
  - K_0.$$
Since $\oo{p_1}$ (respectively $\oo{p_2}$) is a prefix of $\oo{g}$
(respectively $\oo{h}$), it follows easily that $\oo{p_1}\wedge
\oo{p_2} = u$ and $|\oo{p_1}|-|\oo{p_1}\wedge \oo{p_2}| =
d_1(u\pi,p_1)$. Hence
$$\begin{array}{lll}
d_1(u\pi,p_1)&\leq&(d_1(p_1,p_2)+ K_0)N_0 \leq (d_1(p_1,q_j) +
d_1(q_j,q_{j+1}) + d_1(q_{j+1},p_2)+ K_0)N_0\\
&\leq&(2\delta
+1+K_0)N_0
\end{array}$$
and we are done.
\qed

We recall that an automaton is said to be trim if every vertex occurs
in some successful path. Let $\A = (Q,q_0,T,E)$ be a finite trim deterministic
$\widetilde{A}$-automaton recognizing 
$L$ (e.g. the minimal automaton of $L$, see \cite{Ber}). Since $L$ is
factorial, we must have $T = Q$. Let $$\partial
L = \{ \alpha \in \widetilde{A}^{\omega} \mid \alpha^{[n]} \in L \mbox{ for
  every }n \in \N \}.$$
Equivalently, since $\A$ is trim and deterministic, and $T = Q$, we
have
$\partial
L = L_{\omega}(\A)$. Write $\wh{L} = L \cup \partial L$. 
We define a
mapping $d_3: \wh{L} \times \wh{L} \to \mathbb{R}_0^+$ by
$$d_3(\alpha,\beta) = 
\left\{
\begin{array}{ll}
2^{-|\alpha\wedge\beta|}&\mbox{ if } \alpha\neq\beta\\
0&\mbox{ otherwise}
\end{array}
\right.$$
It is immediate that $d_3$ is a distance in $\wh{L}$, indeed an
ultrametric since
$$|\alpha\wedge\gamma| \geq \min\{ |\alpha\wedge\beta|,
|\beta\wedge\gamma| \}$$
holds for all $\alpha,\beta,\gamma \in \wh{L}$. We shall commit a slight
abuse of notation by denoting also by $d_3$ the restriction of $d_3$
to $L \times L$.
 
\bp
\label{buildb}
\bi
\item[(i)] The mutually inverse mappings
$(G,d_2) \to (L,d_3): g \mapsto \oo{g}$ and $(L,d_3) \to
  (G,d_2): u \mapsto u\pi$ are uniformly continuous;
\item[(ii)] $(\wh{L},d_3)$ is the completion of $(L,d_3)$;
\item[(iii)] $(\partial L,d_3)$ is homeomorphic to the boundary of
  $G$.
\ei
\ep

\proof
(i) In view of (\ref{ineq}), it suffices to show that
$$\forall M > 0\; \exists N > 0: \; ((g|h) > N \Rw |\oo{g} \wedge
\oo{h}| > M),$$
$$\forall M > 0\; \exists N > 0: \; (|\oo{g} \wedge
\oo{h}| > N \Rw (g|h) > M).$$ 
Now we apply Lemma \ref{tribound}(iii).

(ii) Let $(\alpha_n)_n$ be a Cauchy sequence in $(\wh{L},d_3)$. For
every $k \in \N$, 
the sequence $(\alpha_n^{[k]})_n$ stabilizes when $n \to
+\infty$. Moreover, $\lim_{n\to +\infty} \alpha_n^{[k]}$ is a prefix
of $\lim_{n\to +\infty} \alpha_n^{[k+1]}$. 
Let $\beta \in A^{\infty}$ be the unique word satisfying $\beta^{[k]}
= \lim_{n\to +\infty} \alpha_n^{[k]}$ for every $k \in \N$. It is
immediate that $\beta \in \wh{L}$ and $\beta = \lim_{n\to +\infty}
\alpha_n$, hence $(\wh{L},d_3)$ is complete. Since $\alpha =
\lim_{n\to +\infty} \alpha^{[n]}$ for every $\alpha \in \partial L$, 
$(\wh{L},d_3)$ is the completion of $(L,d_3)$.

(iii) By (i) and (ii), the uniformly continuous mappings $(G,d_2) \to
(L,d_3): g \mapsto \oo{g}$ and $(L,d_3) \to (G,d_2): u \mapsto u\pi$
admit (unique) continuous extensions to their completions (see
\cite[Section XIV.6]{Dug}), say 
$$\Phi: \wh{G} \to \wh{L}, \quad \Psi: \wh{L} \to \wh{G}.$$
Hence $\Phi\Psi$ is a continuous extension of the identity on $G$ to
its completion $\wh{G}$. Since such an extension is unique, $\Phi\Psi$
must be the identity mapping on $\wh{G}$. Similarly, $\Psi\Phi$
must be the identity mapping on $\wh{L}$ and so $\Phi$ and $\Psi$ are
  mutually inverse homeomorphisms. Therefore the restriction
$\Phi|_{\partial
    G}: \partial G \to \partial L$ must be also a homeomorphism. 
\qed

We have just proved that our construction of $\wh{L}$ constitutes a
model for the hyperbolic completion of $G$. But we must import also to
$\wh{L}$ the algebraic operations of $\wh{G}$ since we shall be
considering homomorphisms soon. Clearly, the binary operation on $L$
is defined as
$$L \times L \to L: (u,v) \mapsto \oo{uv}$$
so that $(G,d_2) \to (L,d_3): g \mapsto \oo{g}$ is also a group
isomorphism. But there is 
another important algebraic operation involved. Indeed, for every $g
\in G$, the left translation $\tau_g:G \to G: x \mapsto gx$ is
uniformly continuous for $d_2$ and so admits a continuous extension
$\wh{\tau}_g: \wh{G} \to \wh{G}$. It follows that the left action of
$G$ in its boundary, $G
\times \partial G \to \partial G: (g,\alpha) \mapsto
\alpha\wh{\tau}_g$, is continuous. We can also replicate this
operation in $\wh{L}$ as follows:

\bp
\label{mix}
Let $u \in L$. Then $\tau_u:L \to L: v \mapsto \oo{uv}$ is
uniformly continuous.
\ep

\proof
It suffices to show that
$$\forall M > 0\; \exists N > 0: \; (|v\wedge w| > N \Rw |\oo{uv} \wedge
\oo{uv}| > M).$$
By Proposition \ref{unibo}(ii), we can take $N = M +  K_0N_0 +
N_0|u|$.
\qed

Therefore $\tau_u$ admits a continuous extension
$\wh{\tau}_u: \wh{L} \to \wh{L}$ and the left action $L
\times \partial L \to \partial L: (u,\alpha) \mapsto
\alpha\wh{\tau}_u$ is continuous. Write $\oo{u\alpha} =
\alpha\wh{\tau}_u$.
For every $\alpha \in
\partial L$, we have
$$\oo{u\alpha} = \oo{u\lim_{n\to +\infty} \alpha^{[n]}} = \lim_{n\to
  +\infty} \oo{u\alpha^{[n]}},$$ hence $(\wh{L},d_3)$ serves as a
model for $(\wh{G},\wh{d_2})$ both topologically and algebraically. From
now on, we shall pursue our work within $(\wh{L},d_3)$.

\section{Uniformly continuous endomorphisms}

We keep all the notation introduced in Section
\ref{sfive}. In particular, $G$ is a finitely generated virtually free
group and $L = \irr\R'$. Following the program announced above, we
work within $(\wh{L},d_3)$.

Given an endomorphism $\p$ of $G$, we denote by $\oo{\p}$ the
corresponding endomorphism of $L$ for the binary operation induced by
the product in $G$, i.e. $u\oo{\p} = \oo{(u\pi)\p}$. To simplify
notation, we shall often write $u\p$ instead of $u\pi\p$ for $u \in
\wt{A}^*$.

We say that $\p$
 satisfies the {\em bounded
  reduction property} if $\{ |u\oo{\p}| - |u\oo{\p} \wedge
(uv)\oo{\p}| : uv \in L \}$ is bounded. In that case, we denote its
maximum by $B_{\p}$. This property was considered originally for free
group automorphisms by Cooper \cite{Coo}.

We fix also the notation $D_{\p} = \max\{ |\oo{a\p}|: a \in \widetilde{A} \}$.

\bt
\label{brlemma}
Let $\p$ be an endomorphism $\p$ of $G$ with finite kernel. Then
$\p$ satisfies the bounded reduction property.
\et

\proof
Suppose that 
$\p$ does not satisfy the bounded
reduction property. Then
$$\forall m \in \N\; \exists u_mv_m \in L :
 |u_m\oo{\p}| - |u_m\oo{\p} \wedge
(u_mv_m)\oo{\p}| > m.$$
Let $X_0 = (K_0 +D_{\p})N_0$. We claim that
\beq
\label{kernel2}
\begin{array}{ll}
\forall m \in \N\; \exists u'_mv'_m \in L :&
( |u'_m\oo{\p}| - |(u'_mv'_m)\oo{\p}| > m\\
&\mbox{and }|(u'_mv'_m)\oo{\p}| -
 |u'_m\oo{\p} \wedge (u'_mv'_m)\oo{\p}| \leq X_0).
\end{array}
\eeq
Indeed, let $m \in \N$. Take $n = m + X_0$ and write
$v_n = a_1\ldots 
a_k$ $(a_i \in \widetilde{A})$. For $i = 0,\ldots,k$, let $w_i =
(u_na_1\ldots a_i)\oo{\p}$. Let $j$ denote the smallest $i$ such that
$|u_n\oo{\p} \wedge w_i| \leq |u_n\oo{\p} \wedge
(u_nv_n)\oo{\p}|$. Take $u'_m = u_n$ and $v'_m = a_1\ldots a_{j-1}$ (since
$j > 0$). Since
$L$ is factorial, we have $u'_mv'_m \in L$. 

Now by minimality of $j$ we get
$$|u_n\oo{\p} \wedge w_{j-1}| > |u_n\oo{\p} \wedge
(u_nv_n)\oo{\p}|.$$
Since $|u_n\oo{\p} \wedge w_j| \leq |u_n\oo{\p} \wedge
(u_nv_n)\oo{\p}|$, it follows that $$|w_{j-1} \wedge w_j| \leq
|u_n\oo{\p} \wedge (u_nv_n)\oo{\p}|.$$ 
Applying Lemma \ref{tribound}(i) to $w_{j-1}\pi$ and $a_j\p$, we get
$$\begin{array}{lll}
|w_{j-1}|&\leq&|w_{j-1} \wedge w_j| + K_0N_0+ N_0|\oo{a_j\p}| \leq 
|w_{j-1} \wedge w_j| + X_0\\
&\leq&|u_n\oo{\p} \wedge (u_nv_n)\oo{\p}| + X_0 < |u_n\oo{\p}|-n + X_0
= |u_n\oo{\p}|-m 
\end{array}$$ 
and so 
$|u'_m\oo{\p}| - |(u'_mv'_m)\oo{\p}| = |u_n\oo{\p}| - |w_{j-1}| > m$.

Suppose that $|w_{j-1}| - |u_n\oo{\p} \wedge w_{j-1}| > X_0$. 
Since we have seen above that $|w_{j-1}| \leq |w_{j-1} \wedge w_j| +
X_0$, we get $|u_n\oo{\p} \wedge w_{j-1}| < |w_{j-1} \wedge w_j|$,
in contradiction with $|w_{j-1} \wedge w_j| \leq |u_n\oo{\p}
\wedge (u_nv_n)\oo{\p}| < |u_n\oo{\p} \wedge w_{j-1}|$. Thus 
$$|(u'_mv'_m)\oo{\p}| -
 |u'_m\oo{\p} \wedge (u'_mv'_m)\oo{\p}| = |w_{j-1}| - |u_n\oo{\p}
 \wedge w_{j-1}| \leq X_0$$ and so (\ref{kernel2}) holds.

% Recalling the notation introduced in (\ref{lote}), we take
% $Y_0 = \max\{ C_{a\oo{\p}} \mid a \in \widetilde{A} \}$. 
We prove that
\beq
\label{kernel3}
\forall m \in \N\; \exists u''_mv''_m \in L :
 |u''_m\oo{\p}| > m \mbox{ and } |(u''_mv''_m)\oo{\p}| \leq X_0 + N_0D_{\p}.
\eeq
Indeed, let $m \in \N$. We have in $\Gamma_A(G)$ geodesics
$$\xymatrix{
1 \ar[rr]^p && g \ar[rr]^q \ar[drr]^r && {u'_m\p} \\
&&&& {(u'_mv'_m)\p}
}$$
where $pq = u'_m\oo{\p}$, $pr = (u'_mv'_m)\oo{\p}$ and $p =
u'_m\oo{\p} \wedge (u'_mv'_m)\oo{\p}$. 
Assume that $u'_m = a_1\ldots a_k$ $(a_i \in \widetilde{A})$.
Let 
$$I = \{ i \in \{ 0,\ldots,k\} \mid \mbox{ there exists a geodesic
  $(a_1\ldots a_{i})\p \mapright{} g \mapright{q} u'_m\p$ in
}\Gamma_A(G) \}.$$  
Clearly, $0 \in I$. We claim that
\beq
\label{kernel4}
(i-1 \in I \mbox{ and } d_1( (a_1\ldots a_{i-1})\p, g) > N_0D_{\p}) \Rw i
\in I
\eeq
holds for $i = 1,\ldots,k$. Indeed, assume that $i-1 \in I$ and
$(a_1\ldots a_{i-1})\p \mapright{y} g \mapright{q} u'_m\p$ is a
geodesic with $y \in L$. Applying Lemma
\ref{bouge}(ii) to the word $a_{i}\inv\oo{\p}$ and the geodesic $yq$, it follows
that there exists some geodesic $(a_1\ldots a_{i-1})\p \mapright{z}
u'_m\p$
such that $z$ and $u$ share a suffix of
length $\geq |yq| - N_0|a_{i}\inv\oo{\p}| \geq |yq|-N_0D_{\p} > |q|$. 
Since $\Gamma_A(G)$ is deterministic, then
our geodesic  
$(a_1\ldots a_{i-1})\p \mapright{z} u'_m\p$ factors through
$g$ and so (\ref{kernel4}) holds.

Since $k \notin I$ due to $|q| > 0$, it follows from (\ref{kernel4})
that $d_1( (a_1\ldots a_{i})\p, g) \leq N_0D_{\p}$ for some $i \in \{
1,\ldots,k\}$.  
let $j$ denote the smallest such $i$. We define
$u''_m = a_{j+1}\ldots a_k$ and $v''_m = v'_m$. Since $L$ is factorial and
$u'_mv'_m \in L$, we have also $u''_mv''_m \in L$. 

By minimality of $j$, we have $d_1( (a_1\ldots a_{i})\p, g) > N_0D_{\p}$ for
$i = 0,\ldots,j-1$. By (\ref{kernel4}), we get $1,\ldots,j \in I$
and so there exists a geodesic
  $(a_1\ldots a_{j})\p \mapright{} g \mapright{q} u'_m\p$ in
$\Gamma_A(G)$. Hence
$$\begin{array}{lll}
|u''_m\oo{\p}|&=&d_1(1,u''_m\p) = d_1((a_1\ldots
a_{j})\p,u'_m\p)\\
&\geq&|q| = |u'_m\oo{\p}| - |(u'_mv'_m)\oo{\p}| > m.
\end{array}$$
Finally,  
$$\begin{array}{lll}
|(u''_mv''_m)\oo{\p}|&=&d_1(1,(u''_mv''_m)\p) = d_1((a_1\ldots
a_{j})\p,(u'_mv'_m)\p)\\
&\leq&d_1((a_1\ldots a_{j})\p,g) + d_1(g,(u'_mv'_m)\p) \leq N_0D_{\p} +
|r|\\
&=&N_0D_{\p} + |(u'_mv'_m)\oo{\p}| -
 |u'_m\oo{\p} \wedge (u'_mv'_m)\oo{\p}| \leq N_0D_{\p} + X_0
\end{array}$$
and so (\ref{kernel3}) holds.

Now, since $|(u''_mv''_m)\oo{\p}|$ is bounded, $u''_mv''_m \in L$ and
$\ker\p$ is finite, then $|u''_mv''_m|$ must be
bounded and so must be  $|u''_m|$. This implies that $|u''_m\oo{\p}|$
must be bounded, contradicting $|u''_m\oo{\p}| > m$. Therefore $\p$
satisfies the bounded reduction property. 
\qed

\bp
\label{uck}
The following conditions are equivalent for a nontrivial endomorphism $\p$ of
$G$:
\bi
\item[(i)] $\p$ is uniformly continuous for $d_2$;
\item[(ii)] {\rm Ker}$\p$ is finite;
\ei
\ep

\proof
(i) $\Rw$ (ii). 
Suppose that $\ker\p$ is infinite. In view of (\ref{ineq}), it
suffices to show that there exists some $\eta > 0$ such that
$$\forall \xi > 0 \; \exists g,h \in G \; (\rho(g,h) < \xi \mbox{ and }
\rho(g\p,h\p) \geq \eta).$$
By (\ref{ineq}), we only need to show that there exists some $M \in
\N$ such that 
$$\forall N \in \N \; \exists g,h \in G \; ((g|h) > N \mbox{ and }
g\p \neq h\p \mbox{ and } ((g\p)|(h\p)) \leq M).$$
Take $M = (g_0\p|1) = 0$ and fix $g_0 \in G \setminus \ker\p$. We
prove the claim by 
showing that 
\beq
\label{kernel1}
\forall N \in \N \; \exists h \in \ker\p :\; ((hg_0)|h) > N.
\eeq
Let $N \in \N$. By Lemma \ref{tribound}(iii), we have
$|\oo{hg_0} \wedge \oo{h}| \leq
((hg_0)|h)$ for every $h \in G$, hence we only need to find out $h \in
\ker\p$ satisfying $|\oo{hg_0} \wedge \oo{h}| > N$. By Lemma
\ref{tribound}(i), we have $|\oo{hg_0} \wedge \oo{h}| \geq |\oo{h}| -
K_0N_0 - N_0|\oo{g_0}|$, 
hence it suffices that $|\oo{h}| > N + K_0N_0+
N_0|\oo{g_0}|$ for some $h \in \ker\p$, and that is ensured by
$\ker\p$ being infinite. Thus 
(\ref{kernel1}) holds as required.

(ii) $\Rw$ (i). Suppose that $\p$ is not uniformly continuous for $d_2$.
In view of (\ref{ineq}), there exists some $\eta > 0$ such that
$$\forall \xi > 0 \; \exists g,h \in G \; (\rho(g,h) < \xi \mbox{ and }
\rho(g\p,h\p) \geq \eta).$$
Hence, by (\ref{ineq}), there exists some $M \in \N$ such that
$$\forall N \in \N \; \exists g,h \in G \; ((g|h) > N \mbox{ and }
g\p \neq h\p \mbox{ and } ((g\p)|(h\p)) \leq M).$$
In view of Lemma \ref{tribound}(iii), we have that
$$
\forall n \in \N \; \exists u_n,v_n \in L \; (|u_n\wedge v_n| > n \mbox{ and }
u_n\oo{\p} \neq v_n\oo{\p} \mbox{ and } |u_n\oo{\p} \wedge v_n\oo{\p}| \leq M).
$$
Let $w_n = u_n\wedge v_n \in L$. Then either $w_n\oo{\p} \neq
u_n\oo{\p}$ or $w_n\oo{\p} \neq v_n\oo{\p}$. Without loss of
generality, we may assume that $w_n\oo{\p} \neq
u_n\oo{\p}$. Suppose that $|w_n\oo{\p}| > M + B_{\p}$. By definition
of $B_{\p}$, we get $|w_n\oo{\p}| - |w_n\oo{\p} \wedge
u_n\oo{\p}| \leq B_{\p}$ and so $|w_n\oo{\p} \wedge
u_n\oo{\p}| > M$. Similarly, $|w_n\oo{\p} \wedge
v_n\oo{\p}| > M$ and so $|u_n\oo{\p} \wedge
v_n\oo{\p}| > M$, a contradiction. Therefore $|w_n\oo{\p}| \leq M +
B_{\p}$ for every $n$. Since $|w_n| > n$ and $L$ is a cross-section
for $\pi$, it follows that $\ker\p$ is infinite.
\qed

Given a uniformly continuous endomorphism $\p$ of $(G,d_2)$, then
$\oo{\p}:L \to L$ is uniformly continuous for $d_3$. Since
$\wh{L}$ is the completion of $(L,d_3)$, then $\oo{\p}$ admits a
unique continuous extension $\Phi:\wh{L} \to \wh{L}$. By continuity,
we have
\beq
\label{exte}
\alpha\Phi = (\lim_{n\to +\infty} \alpha^{[n]})\Phi = \lim_{n\to
  +\infty} \alpha^{[n]}\oo{\p}.
\eeq

\bc
\label{infbrp}
Let $\p$ be a uniformly continuous endomorphism of $G$ and let
$u\alpha \in \partial L$. Then $|u\oo{\p}| - |u\oo{\p} \wedge 
(u\alpha)\Phi| \leq B_{\p}$.
\ec

\proof
We have $(u\alpha)\Phi = \lim_{n\to
  +\infty} (u\alpha^{[n]})\oo{\p}$ by (\ref{exte}). In view of
Proposition \ref{uck}, we have $\lim_{n\to
  +\infty} |(u\alpha^{[n]})\oo{\p}| = +\infty$, hence 
$|u\oo{\p} \wedge (u\alpha)\Phi| = |u\oo{\p} \wedge
(u\alpha^{[m]})\oo{\p}|$ for sufficiently large $m$. Since
$u\alpha^{[m]} \in L$, the claim follows from the definition of $\B_{\p}$.
\qed

\section{Infinite fixed points}

Keeping all the notation and assumptions introduced in the preceding
sections, we fix now a uniformly continuous endomorphism $\p$ of the
finitely generated virtually free group $G$. We adapt
notation introduced in  
\cite{LSil} for free groups, and the proofs are also 
adaptations of proofs in \cite{Sil2}. 

Given
$u \in L$, let $u\sigma = u 
\wedge u\oo{\p}$ and write
$$u = (u\sigma)(u\tau),\quad u\oo{\p} = (u\sigma)(u\rho).$$
Define also
$$u\sigma' = \wedge\{ (uv)\sigma \mid uv \in L \}$$
and write $u\sigma = (u\sigma')(u\sigma'')$.
% Given $u \in \wt{A}^*$ and $n \in \N$, let $us_n$ denote the suffix of
% length $n$ of $u$ if $|u| \geq n$ and $u$ otherwise.

\bl
\label{bsi}
% For every $u \in$ {\rm Pref}$\,${\rm Fix}$\,\Phi$, we have
Let $uv \in L$. Then:
\bi
\item[(i)] $|u\sigma''| \leq B_{\p}$;
\item[(ii)] $|u\sigma| - | u\sigma \wedge (uv)\oo{\p}| \leq
  |u\sigma''|$;
\item[(iii)] $(uv)\oo{\p} =
  (u\sigma')\oo{(u\sigma'')(u\rho)(v\oo{\p})}$;
\item[(iv)] $(uv)\sigma' = (u\sigma')(\displaystyle\bigwedge_{uvz \in L}
  \; (\; \oo{(u\sigma'')(u\rho)((vz)\oo{\p})} \wedge
  (u\sigma'')(u\tau)vz \; )$.
\ei
\el

\proof
(i) We may assume that $|u\sigma| > B_{\p}$.
Let $v$ denote the suffix of length $B_{\p}$ of $u\sigma$ and write
$u\sigma = u'v$. Suppose that $uw \in L$. It suffices to show that
$u'$ is a prefix of $(uw)\oo{\p}$, and this follows from
$$|u'v(u\rho)| - |u'v(u\tau) \wedge (uw)\oo{\p}| = |u\oo{\p}| - |u\oo{\p} \wedge
(uw)\oo{\p}| \leq B_{\p}$$
and $|v| = B_{\p}$.

(ii) Since $u\sigma'$ is a prefix of $u\sigma \wedge (uv)\oo{\p}$.

(iii) Since $u\sigma'$ is a prefix of $(uv)\oo{\p}$ and both sides of
the equality are equivalent in $G$.

(iv) Since $u\sigma'$ is a prefix of $(uv)\sigma'$ by (iii).
\qed

% Let $P_0 = K_0N_0 + V_0 -1$. 
For every $u \in L$, we define
$$u\xi = (% u\sigma' s_{P_0}, 
u\sigma'',u\tau,u\rho,q_0u).$$
Note that there exists precisely one path of the form $q_0 \mapright{u}
q_0u$ in $\A$.

\bl
\label{buiaut}
Let $u,v \in L$ be such that $u\xi = v\xi$ and let $a \in
\widetilde{A}$, $\alpha \in \widetilde{A}^{\infty}$. Then:
\bi
\item[(i)] $ua \in L$ if and only if $va \in L$;
\item[(ii)] if $ua \in L$, then $(ua)\xi = (va)\xi$
\item[(iii)] $\oo{uv\inv} \in \fix\oo{\p}$;
\item[(iv)] $u\alpha \in \wh{L}$ if and only if $v\alpha \in \wh{L}$;
\item[(v)] $u\alpha \in \fix\Phi$ if and only if $v\alpha \in \fix\Phi$;
\item[(vi)] if $\alpha \in \wh{L}$, then $\alpha = \lim_{n\to +\infty}
  \oo{\alpha^{[n]}u}$.
\ei
\el

\proof
(i) Since $u\xi = v\xi$ implies $q_0u = q_0v$.

(ii) Clearly, $q_0u = q_0v$ yields $q_0ua = q_0va$.
Considering $v = a$ in Lemma \ref{bsi}(iii), we may write
$(ua)\sigma = (u\sigma')u'$ and deduce that $u'$, $(ua)\tau$ and
$(ua)\rho$ are all determined by $u\xi$. Hence $(ua)\tau = (va)\tau$,
$(ua)\rho = (va)\rho$ and $u' = v'$.

Finally, since $q_0u = q_0v$, we have $uaz \in L$ if and
only if $vaz \in L$. It follows from Lemma \ref{bsi}(iv) that
there exists a word $x \in L$ which satisfies both $(ua)\sigma' = (u\sigma')x$
and $(va)\sigma' = (v\sigma')x$. Now
$(u\sigma')u' = (ua)\sigma = ((ua)\sigma')((ua)\sigma'') =
(u\sigma')x((ua)\sigma'')$, hence $u' = x((ua)\sigma'')$. Similarly, $v' =
x((va)\sigma'')$. Since $u' = v'$, we get
$(ua)\sigma'' = (va)\sigma''$ and so $(ua)\xi = (va)\xi$. 

(iii) Since
$$\begin{array}{lll}
\oo{(uv\inv)\p}&=&\oo{(u\p)(v\p)\inv} = \oo{(u\sigma)(u\rho)(v\rho)\inv
(v\sigma)\inv} = \oo{(u\sigma)(v\sigma)\inv}\\
&=&\oo{(u\sigma)(u\tau)(v\tau)\inv
(v\sigma)\inv} = \oo{uv\inv}.
\end{array}$$

(iv) We have $u\alpha \in \wh{L}$ if and only if $u\alpha^{[n]} \in L$
for every $n \in \N$. Now we use (i) and induction on $n$.

(v) We have $u\alpha = (u\sigma')(u\sigma'')(u\tau)\alpha$ and in
view of Corollary \ref{infbrp} and (\ref{exte}) also 
$$(u\alpha)\Phi = (u\sigma')\lim_{n\to +\infty}
\oo{(u\sigma'')(u\rho)(\alpha^{[n]}\oo{\p})}.$$ Hence $u\alpha \in
\fix\Phi$ depends just on $u\xi$ and $\alpha$ and we are done.

(vi) Let $m = K_0N_0 + N_0|u|$. By Lemma
\ref{tribound}(i), we have $|\alpha^{[n]} \wedge \oo{\alpha^{[n]}u}|
\geq n-m$ for every $n$, hence 
 $\alpha = \lim_{n\to +\infty}
\alpha^{[n-m]} = \lim_{n\to +\infty}
  \oo{\alpha^{[n]}u}.$
\qed

Given $X \subseteq A^{\infty}$, write
$$\pref X = \{ u \in A^* \mid u\alpha \in X \mbox{ for some } \alpha
\in A^{\infty} \}.$$
Recall the finite trim deterministic
$\widetilde{A}$-automaton $\A = (Q,q_0,Q,E)$ recognizing 
$L$.  We build a (possibly infinite) $\widetilde{A}$-automaton $\A'_{\p} =
(Q',q'_0,T',E')$ by taking 
\bi
\item
$Q' = \{ u\xi \mid u \in \pref\fix\Phi \}$;
\item
$q'_0 = 1\xi$;
\item
$T' = \{ u\xi \in Q' \mid u\tau = u\rho = 1\}$;
\item
$E' = \{ (u\xi,a,v\xi) \in Q' \times \widetilde{A} \times Q' \mid v = ua \in
\pref\fix\Phi \}$.
\ei

Note that $\A'_{\p}$ is deterministic by Lemma \ref{buiaut}(ii) and is
also {\em accessible}: if $u \in \pref\fix\Phi$, then there exists a
path $q'_0 \mapright{u} u\xi$ and so every vertex can be reached from
the initial vertex.

Let 
$S$ denote the set of all vertices $q \in Q'$ such that there exist at
least two edges in $\B'_{\p}$ leaving $q$. Let $Q''$ denote the set of
all vertices $q \in Q'$ such that there exists some path $q
\mapright{} p \in S\cup T'.$ We define $\A''_{\p} =
(Q'',q''_0,T'',E'')$ by taking $q''_0 = q'_0$, $T'' = T' \cap Q''$ and
$E'' = E' \cap (Q'' \times \widetilde{A} \times Q'').$

\bl
\label{sfini}
$S$ is finite.
\el

\proof
In view of Lemma \ref{bsi}, the unique components of $u\xi$ that may
assume infinitely many values are $u\tau$ and $u\rho$. Moreover, we
claim that
\beq
\label{sfini1}
u\tau \neq 1 \Rw |u\rho| \leq B_{\p}
\eeq
holds for every $u \in \pref\fix\Phi$.
Indeed, suppose that $u\tau \neq 1$ and $|u\rho| > B_{\p}$. Write
$\alpha = u\beta$ for some $\alpha \in \fix \Phi$. In view of
Proposition \ref{infbrp}, $|u\rho| > B_{\p}$ yields $|(u\beta)\Phi
\wedge u\oo{\p}| > 
|u\sigma|$ and now $u\tau \neq 1$ yields $((u\beta)\Phi \wedge u\beta)
= (u\oo{\p} \wedge u) = u\sigma$. Since $\beta \neq 1$,
this contradicts $\alpha \in \fix \Phi$. Therefore (\ref{sfini1})
holds.

It is also easy to see that
\beq
\label{sfini2}
|u\rho| > B_{\p} \Rw u\xi \notin S
\eeq
for every $u \in \pref\fix\Phi$. Indeed, if $|u\rho| > B_{\p}$ and $a$
is the first letter of $u\rho$, then, by definition of $B_{\p}$, $(u\sigma)a$ is
a prefix of $(u\alpha)\Phi$ whenever $u\alpha \in \fix\Phi$. Therefore
any edge leaving $u\xi$ in $\A'_{\p}$ must have label $a$ and so
(\ref{sfini2}) holds.

In view of Proposition \ref{uck}, we can define $$W_0 = \max\{ |u|: u
\in L,\; |u\oo{\p}| \leq 2(B_{\p} + D_{\p}-1) \}.$$ 
Let $Z_0 = B_{\p} + N_0(K_0 + W_0)D_{\p}$.
To complete the proof of the lemma, it suffices to prove that
\beq
\label{sfini3}
|u\tau| > Z_0 \Rw u\xi \notin S
\eeq
for every $u \in \pref\fix\Phi$.

Suppose that $|u\tau| > Z_0$ and $(u\xi,a,(ua)\xi), (u\xi,b,(ub)\xi) \in E'$
for some $u \in \pref\fix\Phi$, where $a,b \in \widetilde{A}$ are
distinct. We have $(ua)\xi = v\xi$ for some $v \in \pref\fix\Phi$. By
Lemma \ref{buiaut}(v), we get $ua\alpha \in \fix\Phi$ for some
$\alpha \in \wh{L}$. By (\ref{exte}), we get $ua\alpha = \lim_{n\to
  +\infty} (ua\alpha^{[n]})\oo{\p}$ and so $|(ua\alpha^{[n]})\oo{\p}|
\geq |u|$ for sufficiently large $n$. Let
$$p = \min\{ n \in \N : |(ua\alpha^{[n]})\oo{\p}|
\geq |u| \}.$$
Note that $p > 0$ since $|u\tau| > Z_0$ and by (\ref{sfini1}). Since
$|(ua\alpha^{[p-1]})\oo{\p}| < |u|$ by minimality of $p$, we get 
\beq
\label{sfini4}
|(ua\alpha^{[p]})\oo{\p}| \leq |(ua\alpha^{[p-1]})\oo{\p}| + D_{\p} <
|u| + D_{\p}.
\eeq
On the other hand, 
\beq
\label{sfini5}
|u|-|(ua\alpha^{[p]})\oo{\p} \wedge u| \leq B_{\p},
\eeq
otherwise, by definition of $B_{\p}$, $ua\alpha$ and $(ua\alpha)\Phi$
would differ at position $|(ua\alpha^{[p]})\oo{\p} \wedge u| + 1$.

Similarly, $ub\beta \in \fix\Phi$ for some
$\beta \in \wh{L}$. Defining
$$q = \min\{ n \in \N : |(ub\beta^{[n]})\oo{\p}|
\geq |u| \},$$
we get
\beq
\label{sfini6}
|(ub\beta^{[q]})\oo{\p}| < |u| + D_{\p}
\eeq
and 
\beq
\label{sfini7}
|u|-|(ub\beta^{[q]})\oo{\p} \wedge u| \leq B_{\p}.
\eeq
Write $u = u_1u_2$ with $|u_2| = B_{\p}$. Then by (\ref{sfini4}) and
(\ref{sfini5}) we may write $(ua\alpha^{[p]})\oo{\p} =u_1x$ for some
$x$ such that $|x| < B_{\p} + D_{\p}$. Similarly, (\ref{sfini6}) and
(\ref{sfini7}) yield
$(ub\beta^{[q]})\oo{\p}
= u_1y$ for some
$x$ such that $|x| < B_{\p} + D_{\p}$. Writing $w =
\oo{(\beta^{[q]})\inv b\inv a\alpha^{[p]}}$, it follows that
$w\p = (y\inv x)\pi$ and so $|w\oo{\p}| \leq 2(B_{\p} +
D_{\p}-1)$. Hence $|w| \leq W_0$. Applying Lemma \ref{tribound}(i) to
$g = (ub\beta^{[q]})\pi$ and $h = w\pi$, we get
$$|ub\beta^{[q]}| \leq |ub\beta^{[q]} \wedge ua\alpha^{[p]}| +
N_0(K_0 + |w|) \leq |u| + N_0(K_0 + W_0)$$
and so $q < N_0(K_0 + W_0)$. Hence, in view of (\ref{sfini1}), we get
$$\begin{array}{lll}
|u\tau|&=&|u|-|u\sigma| \leq |(ub\beta^{[q]})\oo{\p}|-|u\sigma| \leq
|u\oo{\p}| + |(b\beta^{[q]})\oo{\p}|-|u\sigma|\\
&\leq&|u\rho| + N_0(K_0 + W_0)D_{\p} \leq B_{\p} + N_0(K_0 + W_0)D_{\p},
\end{array}$$
contradicting
$|u\tau| > Z_0$. Thus (\ref{sfini3}) holds and the lemma is proved.
\qed

We say that an infinite fixed point $\alpha \in \fix\Phi \cap \partial
L$ is {\em singular} if $\alpha$ belongs to the topological closure
$(\fix\p)^c$ of $\fix\p$. Otherwise, $\alpha$ is said to be {\em regular}.
We denote by $\sing\Phi$ (respectively $\reg\Phi$) the set of all
singular (respectively regular) infinite fixed points of $\Phi$.

\bt
\label{afini}
Let $\p$ be a uniformly continuous endomorphism of a finitely
generated virtually free
group $G$. Then:
\bi
\item[(i)] the automaton $\A''_{\p}$ is finite;
\item[(ii)] $L(\A''_{\p}) =$ {\rm Fix}$\,\oo{\p}$;
\item[(iii)] $L_{\omega}(\A''_{\p}) =$ {\rm Sing}$\,\Phi$.
\ei
\et

\proof
(i) The set $T'$ is finite and $S$ is finite by Lemma \ref{sfini}. On
the other hand, by definition of $S$, there are only finitely many
paths in $\A'_{\p}$ of the form $p \mapright{} q$ with $p,q \in S \cup
T' \cup \{ q'_0\}$ and no intermediate vertex in $S \cup
T' \cup \{ q'_0\}$. Therefore $Q''$ is finite and so is $\A''_{\p}$.

(ii) Every $u \in L$ labels at most a unique path $q'_0 = 1\xi
\mapright{u} u\xi$ out 
of the initial vertex in $\A'_{\p}$. On the other hand, if $q'_0 = 1\xi
\mapright{u} q'$ is a path in $\A'_{\p}$, then the fourth component of
$\xi$ yields a path $q_0 \mapright{u} q$ in $\A$ and so $u \in
L$. Hence 
$$L(\A'_{\p}) = \{ u \in L \mid u\xi \in T'\} = \{ u \in L \mid u\tau
= u\rho = 1\} = \fix\oo{\p}.$$
Since $L(\A''_{\p}) = L(\A'_{\p})$, (ii) holds.

(iii) Let $\alpha \in L_{\omega}(\A''_{\p})$. Then there exists some
$q'' \in Q''$ and some infinite sequence $(i_n)_n$ such that $q''_0
\vlongmapright{\alpha^{[i_n]}} q''$ is a path in $\A''_{\p}$ for every
  $n$. Write $u = \alpha^{[i_1]}$ and let $v_n =
  \oo{\alpha^{[i_n]}u\inv}$. By Lemma \ref{buiaut}(iii), we have
$v_n \in \fix\oo{\p}$ for every $n$. It follows from Lemma
\ref{buiaut}(vi) that $\alpha = \lim_{n\to +\infty} v_n$, thus
$\alpha \in \sing\Phi$. 

Conversely, let $\alpha \in \sing\Phi$. Then we may write $\alpha =
\lim_{n\to +\infty} v_n$ for some sequence $(v_n)_n$ in
$\fix\oo{\p}$. Let $k \in \N$. For large enough $n$, we have
$\alpha^{[k]} = v_n^{[k]}$ and so there is some path $$q''_0
\longmapright{\alpha^{[k]}} q''_k \mapright{w} t''_k \in T'',$$
where $\alpha^{[k]}w = v_n$. Thus $\alpha \in L_{\omega}(\A''_{\p})$
as required.
\qed

Recall now the continuous extensions $\wh{\tau}_u: \wh{L} \to \wh{L}$
of the uniformly
continuous mappings $\tau_u:L \to L: v \mapsto \oo{uv}$ defined
for each $u \in L$ (see Proposition \ref{mix}). As remarked before,
this is equivalent to say that the left action $L
\times \partial L \to \partial L: (u,\alpha) \mapsto
\oo{u\alpha}$ is continuous. Identifying $L$ with $G$ and
$\partial L$ with $\partial G$, we have a continuous action (on the
left) of $G$ on $\partial G$. Clearly, this action restricts to a left
action of $\fix\p$ on $\fix\Phi \cap \partial G$: if $g\in \fix\p$ and
$\alpha \in \fix\Phi \cap \partial G$, with $\alpha = \lim_{n\to
  +\infty} g_n$ $(g_n \in G)$, then
$$\begin{array}{lll}
(g\alpha)\Phi&=&(g \lim_{n\to
  +\infty} g_n)\Phi = (\lim_{n\to
  +\infty} gg_n)\Phi = \lim_{n\to
  +\infty} (gg_n)\p\\
&=&\lim_{n\to
  +\infty} (g\p)(g_n\p) = (g\p)\lim_{n\to
  +\infty} g_n\p = g(\lim_{n\to  +\infty} g_n)\Phi\\
&=&g(\alpha\Phi) =
g\alpha.
\end{array}$$
Moreover, the $(\fix\p)$-orbits of $\sing\Phi$ and $\reg\Phi$ are
disjoint: if $\alpha \in \sing\Phi$, we can write $\alpha = \lim_{n\to
  +\infty} g_n$ with the $g_n \in \fix\p$ and get $g\alpha = \lim_{n\to
  +\infty} gg_n$ with $gg_n \in \fix\p$ for every $n$; hence $\alpha
\in \sing\Phi \Rw g\alpha
\in \sing\Phi$ and the action of $g\inv$ yields the converse implication.
 
We can now prove the main result of this section.

\bt
\label{groupifp}
Let $\p$ be a uniformly continuous endomorphism of a finitely
generated virtually free
group $G$. Then {\rm Reg}$\,\Phi$ has finitely many ({\rm Fix}$\,\p)$-orbits.
\et

\proof
Let $P$ be the set of all infinite paths $s'_0 \mapright{a_1} s'_1
\mapright{a_2} \ldots$ in
$\A'_{\p}$ such that:
\bi
\item
$s'_0  \in S \cup \{ q_0\}$;
\item
$s'_n \notin S \cup \{ q_0\}$ for every $n > 0$;
\item
$s'_n \neq s'_m$ whenever $n \neq m$.
\ei
By Lemma \ref{sfini}, there are only
finitely many choices for $s'_0$. Since $A$ is finite and $\A'_{\p}$
is deterministic, there are only
finitely many choices for $s'_1$, and from that vertex onwards, the
path is univocally determined due to $s'_n \notin S$ $(n\geq 1)$. Hence $P$ is
finite, and we may assume that it consists of paths $p'_i
\mapright{\alpha_i} \ldots$ for $i = 1,\ldots,m$. Fix a path $q'_0
\mapright{u_i} p_i$ for each $i$ and let $X = \{ u_1\alpha_1, \ldots,
u_m\alpha_m \} \subseteq \partial L$. We claim that $X \subseteq
\reg\Phi$. 

Let $i \in \{ 1,\ldots,m \}$ and write $\beta = u_i\alpha_i$. To show
that $\beta \in \fix\Phi$, it 
suffices to show that $\lim_{n \to +\infty} \beta^{[n]}\oo{\p} =
\beta$. Let $k \in \N$. We must show that there exists some $r \in
\N$ such that
\beq
\label{groupifp1} 
n \geq r \Rw |\beta^{[n]}\oo{\p} \wedge \beta| > k.
\eeq
In view of Proposition \ref{uck}, there exists some $r > k$ such that 
$$n \geq r \Rw |\beta^{[n]}\oo{\p}| > k + B_{\p}.$$
Suppose that $|\beta^{[n]}\oo{\p} \wedge \beta| \leq k$ for some $n
\geq r$. Then $|\beta^{[n]}\sigma| \leq k$. Since $k < r \leq n$, it
follows that $\beta^{[n]}\tau \neq 1$. On the other hand, since
$|\beta^{[n]}\oo{\p}| > k + B_{\p}$, we get $|\beta^{[n]}\rho| >
B_{\p}$. In view of (\ref{sfini1}), this contradicts $\beta^{[n]}\xi
\in Q'$. Therefore (\ref{groupifp1}) holds for our choice of $r$ and
so $X \subseteq
\fix\Phi$. Since the path $q'_0 \mapright{\beta} \ldots$ can visit only
finitely often a given vertex, then $\beta \notin
L_{\omega}(\A''_{\p})$ and so $X \subseteq
\reg\Phi$ by Theorem \ref{afini}(iii). 

By the previous comments on $(\fix\p)$-orbits, the $(\fix\p)$-orbits
of the elements of $X$ must be contained in $\reg\Phi$.
We complete the proof of the theorem by proving the opposite
inclusion.

Let $\beta \in \reg\Phi$. By Theorem \ref{afini}(iii), we have
$\beta \notin L_{\omega}(\A''_{\p})$ and so there exists a
factorization $\beta =u\alpha$ and a path 
$$q'_0 \mapright{u} p' \mapright{\alpha} \ldots$$
in $\A'_{\p}$ such that $p'$ signals the last occurrence of a vertex from $S
\cup \{ q'_0\}$. We claim that no vertex is repeated after
$p'$. Otherwise, since no vertex of $S$ appears after $p'$, we would
get a factorization of $p' \mapright{\alpha} \ldots$ as
$$p' \mapright{v} q' \mapright{w} q' \mapright{w} \ldots$$
and by Lemma \ref{buiaut}(iii) and (iv) we would get
$(uvw^nv\inv u\inv)\pi \in \fix\p$ and $$\beta = \lim_{n \to +\infty}
\oo{uvw^nv\inv u\inv},$$
contradicting $\beta \in \reg\Phi$. Thus no
vertex is repeated after $p'$ and so we must have $p' = p'_i$ and
$\alpha = \alpha_i$ for some $i \in \{ 1,\ldots, m\}$. It follows
that $\beta = u\alpha_i$. By Lemma \ref{buiaut}(iii), we get
$\oo{uu_i\inv} \in \fix\oo{\p}$ and we are done.
\qed

Theorem \ref{groupifp} is somehow a version for infinite fixed points of
Theorem \ref{gtholds}, which we proved before for finite fixed
points. Note however that $\sing\Phi$ has {\em not} in general
finitely many $(\fix\p)$-orbits since  $\sing\Phi$ may be uncountable
(take for instance the identity automorphism on a free group of rank 2).

Since every finite set is closed in a metric space, we obtain the
following corollary from Theorem \ref{groupifp}:

\bc
\label{finfix}
Let $\p$ be a uniformly continuous endomorphism of a finitely
generated virtually free
group $G$ with $\fix\p$ finite. Then {\rm Fix}$\,\Phi$  is finite.
\ec

\section{Classification of the infinite fixed points}

We can now investigate the nature of the infinite fixed points of
$\Phi$ when $\p$ is an automorphism. Since both $\p$ and $\p\inv$ are
then uniformly continuous by Proposition \ref{uck}, they extend to
continuous mappings $\Phi$ and $\Psi$ which turn out to be mutually
inverse in view of the uniqueness of continuous extensions to the
completion. Therefore
$\Phi$ is a bijection.
% The {\em attraction basin} of $\alpha$ is
% $$\att(\alpha) =  \{ \beta \in G \cup \partial_AG \mid \alpha \in
% \ad(\beta\Phi^n)_n \}.$$
We say that $\alpha \in \reg\Phi$ is:
\begin{itemize}
\item
an {\em attractor} if 
$$\exists \varepsilon > 0 \; \forall \beta \in \wh{L}\;
(d_3(\alpha,\beta) < \varepsilon \Rw \lim_{n\to +\infty} \beta\Phi^n =
\alpha).$$
\item
a {\em repeller} if 
$$\exists \varepsilon > 0 \; \forall \beta \in \wh{L}\;
(d_3(\alpha,\beta) < \varepsilon \Rw \lim_{n\to +\infty} \beta\Phi^{-n} =
\alpha).$$
\end{itemize}
The latter amounts to say that $\alpha$ is an attractor for
$\Phi\inv$.
There exist other types but they do not occur in our context as we
shall see.

We say that an attractor $\alpha \in
\reg\Phi$ is {\em exponentially stable} if
$$\exists \varepsilon,k,\ell > 0 \; \forall \beta \in \wh{L}\;
\forall n \in \N \;
(d_3(\alpha,\beta) < \varepsilon \Rw d_3(\alpha,\beta\Phi^n) \leq
k2^{-\ell n}d_3(\alpha,\beta)).$$
This is equivalent to say that
\beq
\label{exps}
\exists M,N,\ell > 0 \; \forall \beta \in \wh{L}\;
\forall n \in \N \;
(|\alpha \wedge \beta| > M \Rw |\alpha \wedge \beta\Phi^n| + N > \ell
n + |\alpha \wedge \beta|).
\eeq

A repeller $\alpha \in
\reg\Phi$ is exponentially stable if it is an exponentially stable
attractor for
$\Phi\inv$.

\bt
\label{classi}
Let $\p$ be an automorphism of a finitely
generated virtually free
group $G$. Then {\rm Reg}$\,\Phi$ contains only exponentially stable
attractors and exponentially stable repellers.
\et

\proof
Let $\alpha \in
\reg\Phi$ and write $\alpha = a_1a_2\ldots$ with $a_i \in
\widetilde{A}$. Then there exists a path
$$1\xi \mapright{a_1} \alpha^{[1]}\xi \mapright{a_2} \alpha^{[2]}\xi
\mapright{a_3} \ldots$$ 
in $\A'_{\p}$. Let $Y_0 = B_{\p}(D_{\p\inv}+1) + B_{\p\inv}(D_{\p}+1)$ and let 
$$V = \{ u\xi \in Q': |u\tau| > Y_0 \mbox{ or } |u\rho| > Y_0 \}.$$
It is easy to see that $Q' \setminus V$ is finite. We saw in the proof
of Theorem \ref{groupifp} that there are only finitely many
repetitions of vertices in a path in $\A'_{\p}$ labelled by a regular
fixed point, hence there exists some $n_0 \in \N$ such that
\beq
\label{last}
\alpha^{[n]}\xi \in V \mbox{ for every }n \geq n_0.
\eeq
Now we consider two
cases:

\medskip

\noindent
\underline{Case I}: $\alpha^{[n_0]}\tau = 1$. 

\medskip

We claim that 
\beq
\label{classi5}
\alpha^{[n]}\tau = 1 \mbox{ for every } n \geq n_0.
\eeq
The case $n = n_0$ holds in Case I, so assume that $\alpha^{[n]}\tau =
1$ for some $n \geq n_0$. Then $\alpha^{[n]} \in V$ and so
$|\alpha^{[n]}\rho| > Y_0 > 2B_{\p}$. Since $|\alpha^{[n+1]}\oo{\p}|
\geq |\alpha^{[n]}\oo{\p}| - B_{\p}$ by definition
of $B_{\p}$, then 
% $\alpha^{[n]}$ is a prefix
% of $\alpha^{[n+1]}\oo{\p}$ and so $|\alpha^{[n+1]}\tau| \leq
% 1$. Moreover, 
$$\begin{array}{lll}
|\alpha^{[n+1]}\rho|&\geq&|\alpha^{[n+1]}\oo{\p}| - |\alpha^{[n+1]}|
\geq |\alpha^{[n]}\oo{\p}| - B_{\p} - |\alpha^{[n]}| -1 =
|\alpha^{[n]}\rho| - B_{\p} -1\\
&>&Y_0 - B_{\p}
-1 > B_{\p}.
\end{array}$$
By (\ref{sfini1}), we get $\alpha^{[n+1]}\tau = 1$ and so
(\ref{classi5}) holds.

Next we show that
\beq
\label{classi6}
((\alpha^{[n]}\gamma)\Phi)^{[n+1]} = \alpha^{[n+1]}
\eeq
if $n\geq n_0$ and $\alpha^{[n]}\gamma \in \wh{L}$.
Indeed, by (\ref{classi5}) we have
$\alpha^{[n]}\oo{\p} = \alpha^{[n]}(\alpha^{[n]}\rho)$ and
$|\alpha^{[n]}\rho| > Y_0 > B_{\p}$. By the definition
of $B_{\p}$ and Corollary \ref{infbrp}, we get
$((\alpha^{[n]}\gamma)\Phi)^{[n+1]} =
\alpha^{[n]}(\alpha^{[n]}\rho)^{[1]}$. Considering the particular case
$\gamma = a_{n+1}$, we also get««
$$(\alpha^{[n+1]}\oo{\p})^{[n+1]} =
\alpha^{[n]}(\alpha^{[n]}\rho)^{[1]} =
((\alpha^{[n]}\gamma)\Phi)^{[n+1]}.$$
Since $\alpha^{[n+1]}\tau = 1$ by (\ref{classi5}), we have
$(\alpha^{[n+1]}\oo{\p})^{[n+1]} = \alpha^{[n+1]}$ and so (\ref{classi6}) holds.

Hence we may write $(\alpha^{[n]}\gamma)\Phi = \alpha^{[n+1]}\gamma'$
whenever $\alpha^{[n]}\gamma \in \wh{L}$. Iterating, it follows that,
for all $k \geq n_0$ and $n \in \N$, 
$\alpha^{[k]}\gamma \in \wh{L}$ implies
$(\alpha^{[k]}\gamma)\Phi^n = \alpha^{[k+n]}\gamma'$ for some 
$\gamma'$. By considering $\beta = \alpha^{[k]}\gamma$ and
$\alpha^{[k]} = \alpha \wedge \beta$, we deduce that
$$|\alpha \wedge \beta| \geq n_0 \Rw |\alpha \wedge \beta\Phi^n| \geq
n + |\alpha \wedge \beta|$$ holds for all $\beta \in \wh{L}$ and $n
\in \N$. Therefore (\ref{exps}) holds and so $\alpha$ is an exponentially
stable attractor in this case.

\medskip

Now, if $|\alpha^{[t]}\tau| = 1$ for some $t > n_0$, we can always
replace $n_0$ by $t$ and deduce by Case I that $\alpha$ is an exponentially
stable attractor. Thus we may assume that:

\medskip

\noindent
\underline{Case II}: $\alpha^{[n]}\tau \neq 1$ for every $n \geq n_0$.

\medskip

By replacing $n_0$ by a larger integer if necessary, we may assume
that (\ref{last}) is also satisfied when we consider the equivalents
of $\xi$ and $V$ for $\p\inv$.

Since $\p$ is injective, there exists some $n_1 \geq n_0$ such that
$|\alpha^{[n_1]}\oo{\p}| \geq n_0 + B_{\p}$. 
Since $\alpha^{[n_1]}\tau \neq 1$, it follows from (\ref{sfini1}) that
$|\alpha^{[n_1]}\rho| \leq B_{\p}$, hence $\alpha^{[n_1]}\sigma =
\alpha^{[n_2]}$ for some $n_2 \geq n_0$. Write $x =
\alpha^{[n_1]}\rho$. Then $\alpha^{[n_1]}\oo{\p} = \alpha^{[n_2]}x$
yields $\alpha^{[n_1]} =
\oo{(\alpha^{[n_2]}\oo{\p\inv})(x\oo{\p\inv})}$ and so
$$n_1 = |\alpha^{[n_1]}| \leq |\alpha^{[n_2]}\oo{\p\inv}| +
|x\oo{\p\inv}| \leq |\alpha^{[n_2]}\oo{\p\inv}| + B_{\p}D_{\p\inv}.$$
% < |\alpha^{[n_2]}\oo{\p\inv}| + Y_0.$$
On the other hand, $|\alpha^{[n_1]}\rho| \leq B_{\p} < Y_0$ and
$\alpha^{[n_1]} \in V$ together yield $Y_0 < |\alpha^{[n_1]}\tau| =
n_1 - n_2$ and so
$$n_2 + B_{\p\inv} < n_1 - Y_0 + B_{\p\inv} < n_1 - B_{\p}D_{\p\inv}
\leq |\alpha^{[n_2]}\oo{\p\inv}|.$$ 
In view of (\ref{sfini1}), we can apply Case I to $\p\inv$, hence
$\alpha$ is an exponentially 
stable attractor for $\p\inv$ and therefore an exponentially
stable repeller for $\p$.
\qed

\section{Example and open problems}

We include a simple example which illustrates some of the
constructions introduced earlier:

\bigskip

\noindent
{\bf Example.} Let $G = \Z \times \Z_2$ and let $A = \{ a,b,c\}$. Note
that this is not the canonical set of generators, which would not
work. Then
the matched homomorphism $\pi:\wt{A}^* \to G$ defined by
$$a\pi = (1,0),\quad b\pi = (0,1), \quad c\pi = (1,1)$$
yields
$$\geo_A(G) = (a \cup c)^* \cup (a\inv \cup c\inv)^* \cup \{ b, b\inv\}$$
and we can take 
$$\begin{array}{lll}
\R&=&\{ (xx\inv,1) \mid x \in \wt{A} \} \cup \{
(a^{\varepsilon}b^{\delta},c^{\varepsilon}),
(b^{\delta}a^{\varepsilon},c^{\varepsilon}),
(c^{\varepsilon}b^{\delta},a^{\varepsilon}), 
(b^{\delta}c^{\varepsilon},a^{\varepsilon})  \mid \delta, \varepsilon =
\pm 1 \}\\
&\cup&\{ (ac\inv,b), (c\inv a,b), (a\inv c,b), (ca\inv,b), (b^2,1),
(b^{-2},1) \}
\end{array}$$
to get $\geo_A(G) = \irr \R$. Ordering $\wt{A}$ by $a < c < a\inv <
c\inv < b < b\inv$, we get
$$L = a^*(1 \cup c) \cup (a\inv)^*(1 \cup c\inv) \cup b,$$
recognized by the automaton $\A$ depicted by
$$\xymatrix{
& q_0 \ar@{<->}[l] \ar[rrr]^a \ar[drrr]^{b,c,c\inv} \ar[d]_{a\inv} &&& q_1
\ar[d]^c \ar[r] \ar@(ul,ur)^a& \\
& q_2 \ar[l] \ar[rrr]_{c\inv} \ar@(dl,dr)_{a\inv} &&& q_3 \ar[r] &
}$$
Hence $\partial L = L_{\omega}(\A) = \{ a^{\omega}, (a\inv)^{\omega}
\}$.

Let $\p$ be the endomorphism of $G$ defined by $(m,n)\p =
(2m,n)$. Then $\p$ is injective and therefore uniformly continuous,
admiting a continuous extension $\Phi$ to $\wh{L}$. Since $B_{\p} =
0$, it is easy to check that $\A'_{\p}$ is the automaton 
$$\xymatrix{
&&& b\xi \ar[r] &&& \\
\ldots & a^{-2}\xi \ar[l]^{a\inv} & a^{-1}\xi \ar[l]^{a\inv} & 1\xi
\ar[l]^{a\inv} \ar[r]_a \ar[u]_b \ar@{<->}[ul] & a\xi \ar[r]_a &
a^2\xi \ar[r]_a & \ldots
}$$
and $$1\xi = (1,1,1,q_0), \quad b\xi = (1,1,1,q_3), \quad a^n\xi =
(1,1,a^n,q_1), \quad a^{-n}\xi = (1,1,a^{-n},q_2)$$
for $n \geq 1$. Note that in general we ignore how to compute
$\A'_{\p}$, our proofs being far from constructive!

It is immediate that $\fix\Phi = \{ 1,b, a^{\omega}, (a\inv)^{\omega}
\}$. Moreover, the regular infinite fixed points $a^{\omega}$ and
$(a\inv)^{\omega}$ are both exponentially
stable attractors.

\bigskip

Finally, we end the paper with some easily predictable open problems:

\bq
Is it possible to generalize Theorems \ref{gtholds}, \ref{groupifp}
and \ref{classi} to 
arbitrary finitely generated hyperbolic groups? 
\eq

Paulin proved that Theorem \ref{gtholds} holds for automorphisms of
hyperbolic groups \cite{Pau}.

\bq
Is $\fix\p$ effectively computable when $\p$ is an endomorphism of a
finitely generated virtually free group?
\eq

For the moment, only the case of free group automorphisms is known
(Maslakova, \cite{Mas}).

\section*{Acknowledgements}

The author acknowledges support from the European Regional Development
Fund through the 
programme COMPETE and from the Portuguese Government through FCT --
Funda\c c\~ao para a Ci\^encia e a Tecnologia, under the project
PEst-C/MAT/UI0144/2011.

\end{document}